\documentclass{birkjour}

\usepackage{lscape,latexsym,amssymb}
\usepackage{frenchineq}
\usepackage{pgf}
\usepackage{tikz}
\usetikzlibrary{calc}
\usetikzlibrary{arrows}
\usetikzlibrary{shapes} %
\usetikzlibrary{shapes}
\usetikzlibrary{plotmarks}
\usepackage{accents}
\usepackage{stackengine}

\usepackage[utf8]{inputenc}
\usepackage[english]{babel}
\usepackage[T1]{fontenc}
\usepackage{calc,bm}
\usepackage{tabularx,enumitem,natbib,nameref,hyperref}
\usepackage[mathscr]{eucal}
\usepackage{cleveref}
\usepackage[colorinlistoftodos,bordercolor=orange,backgroundcolor=orange!20,linecolor=orange,textsize=scriptsize]{todonotes}
\usepackage{mathtools}
\DeclarePairedDelimiter\ceil{\lceil}{\rceil}
\DeclarePairedDelimiter\floor{\lfloor}{\rfloor}
\renewcommand{\b}{\mathbf}

\newcommand{\NN}{\mathbb{N}}

\newcommand{\RR}{\mathbb{R}}

\newcommand{\mR}{\mathbb{R}_\bot}
\newcommand{\pR}{\mathbb{R}_\top}
\newcommand{\pmR}{\overline{\mathbb{R}}}
\newcommand{\pmRX}{\pmR^\X}
\newcommand{\pmRXmax}{\pmR^\X_{\max}}
\newcommand{\pmRXmin}{\pmR^\X_{\min}}
\newcommand{\pmRmax}{\pmR_{\max}}
\newcommand{\pmRmin}{\pmR_{\min}}

\newcommand{\G}{\mathcal{G}}

\newcommand{\bB}{\bar B}
\newcommand{\bF}{\bar F}
\newcommand{\cL}{\mathcal{L}}
\newcommand{\cI}{\mathcal{I}}

\newcommand{\Hk}{\mathscr{H}_k}

\newcommand{\Hsc}{\mathscr{H}}
\newcommand{\Asc}{\mathscr A}

\makeatletter
\newcommand{\dotminus}{\mathbin{\text{\@dotminus}}}

\newcommand{\@dotminus}{%
	\ooalign{\hidewidth\raise1ex\hbox{.}\hidewidth\cr$\m@th-$\cr}%
}
\newcommand{\dotbminus}{\mathbin{\text{\@dotbminus}}}

\newcommand{\@dotbminus}{%
	\ooalign{\hidewidth\raise-0.2ex\hbox{.}\hidewidth\cr$\m@th-$\cr}%
}

\newcommand{\dotbplus}{\mathbin{\text{\@dotbplus}}}

\newcommand{\@dotbplus}{%
	\ooalign{\hidewidth\raise-0.6ex\hbox{.}\hidewidth\cr$\m@th+$\cr}%
}

\makeatother

\newcommand{\X}{\mathscr{X}}

\newcommand{\Z}{\mathscr{Z}}

\newcommand{\Cop}{C^{op}}

\newcommand{\bc}{b_{\operatorname{conv}}}%
\newcommand{\bsc}{b_{\operatorname{sconv}}}%
\newcommand{\bL}{b_{\operatorname{lip}}}%

\DeclareMathOperator*{\Sp}{span}
\DeclareMathOperator*{\Rg}{Rg}

\DeclareMathOperator*{\Lipb}{Lip_c}

\DeclareMathOperator*{\Id}{Id}
\newcommand{\F}{\mathscr{F}} %

\renewcommand{\d}{\mathrm{d}} %
\newcommand{\tb}{\textbf} %
\newcommand{\R}{\RR} %
\newcommand{\N}{\NN} %
\renewcommand{\b}{\mathbf}

\newcommand{\wrT}{\text{w.r.t.}}

\theoremstyle{plain}
\newtheorem{Theorem}{Theorem}[section]
\newtheorem*{Theorem*}{Theorem}
\newtheorem{Corollary}[Theorem]{Corollary}
\newtheorem{Proposition}[Theorem]{Proposition}

\newtheorem{Lemma}[Theorem]{Lemma}

\theoremstyle{definition}
\newtheorem{Definition}[Theorem]{Definition}
\newtheorem{Example}[Theorem]{Example}%
\newtheorem{Remark}[Theorem]{Remark}%
\date{\today}

\makeatletter
\let\orgdescriptionlabel\descriptionlabel
\renewcommand*{\descriptionlabel}[1]{%
	\let\orglabel\label
	\let\label\@gobble
	\phantomsection
	\edef\@currentlabel{#1\unskip}%
	\let\label\orglabel
	\orgdescriptionlabel{#1}%
}
\makeatother

\catcode`\@=11
\renewcommand{\subjclassname}{%
	\textup{1991} Mathematics Subject Classification}
\@xp\let\csname subjclassname@1991\endcsname \subjclassname
\@namedef{subjclassname@2000}{%
	\textup{2000} Mathematics Subject Classification}
\@namedef{subjclassname@2010}{%
	\textup{2010} Mathematics Subject Classification}
\@namedef{subjclassname@2020}{%
	\textup{2020} Mathematics Subject Classification}
\def\@settitle{\begin{center}%
		\baselineskip14\p@\relax
		\bfseries
		\uppercasenonmath\@title
		\@title
	\end{center}%
}
\catcode`\@=12
\title{Tropical reproducing kernels and optimization}
\author[Aubin-Frankowski]{Pierre-Cyril Aubin-Frankowski}
\address{INRIA and Département d’Informatique, École Normale Supérieure, 
	PSL Research University}
\email{pierre-cyril.aubin@inria.fr}
\author[Gaubert]{Stéphane Gaubert}%

\address{INRIA and CMAP, \'Ecole polytechnique, IP Paris, CNRS}
\email{stephane.gaubert@inria.fr}
\makeatletter
\@namedef{subjclassname}{%
	\textup{2020} Mathematics Subject Classification}
\makeatother

\date{\today}%

\begin{document}
	\maketitle
	\begin{abstract}
		Hilbertian kernel methods and their positive semidefinite kernels have been extensively used in various fields of applied mathematics and machine learning, owing to their several equivalent characterizations.
		We here unveil an analogy with concepts from tropical geometry, proving that tropical positive semidefinite kernels are also endowed with equivalent viewpoints, stemming from Fenchel-Moreau conjugations. This tropical analogue of Aronszajn's theorem shows that these kernels correspond to a feature map, define monotonous operators, and generate max-plus function spaces endowed with a reproducing property. They furthermore include all the Hilbertian kernels classically studied as well as Monge arrays. However, two relevant notions of tropical reproducing kernels must be distinguished, based either on linear or sesquilinear interpretations. The sesquilinear interpretation is the most expressive one, since reproducing spaces then encompass classical max-plus spaces, such as those of (semi)convex functions. In contrast, in the linear interpretation, the reproducing kernels are characterized by a restrictive condition, von Neumann regularity. Finally, we provide a tropical analogue of the ``representer theorems'', showing that a class of infinite dimensional regression and interpolation problems admit solutions lying in finite dimensional spaces. We illustrate this theorem by an application to optimal control, in which tropical kernels allow one to represent the value function.
		
	\end{abstract}	
	
	\section{Introduction}
	
	\paragraph{\bf Context.} Since the foundation of their theory \citep[see][for a historical summary]{aronszajn50theory}, reproducing kernel Hilbert spaces (RKHSs) have played an eminent role among linear function spaces \citep{saitoh16theory}, all the more in machine learning \citep{steinwart08support, scholkopf02learning}. However, optimization problems frequently involve some more intricate function spaces, requiring dedicated structures \citep{Pallaschke1997}, such as the space of convex functions used in convex regression \citep[e.g.][]{Seijo2011}. For instance, value functions, the solutions of the Hamilton-Jacobi-Bellman (HJB) equation in optimal control, are generically nonsmooth while still being semiconcave \citep{Cannarsa2004}. Moreover, Lax-Oleinik semigroups, i.e., evolution semigroups of HJB equations with a Hamiltonian convex in the adjoint variable, are tropically linear. This fact spurred research in tropical / idempotent functional analysis \citep{KM,McEneaney2006}. At the price of a change to $(\max,+)$ operations, many concepts, e.g.\ operators, have been defined by analogy with the linear setting. This allowed to develop new classes of numerical methods for the HJB equation~\citep{a5,asma-art-1,McEneaney2006,curseofdim,Dower2015}.
	Furthermore, several research directions have been explored at the interface between tropical geometry, probablity theory and machine learning. These include studies of the tropicalization of stochastic processes~\citep{bellman} or
	of Gaussian measures~\citep{Tran2020}, tropical support vector machines~\citep{yoshida2021TropicalSV}, tropical principal component analysis~\citep{yoshidatropPCAphylo}, quantification of the expressivity of deep neural networks~\citep{pmlr-v80-zhang18i,2104.08135} or their approximation~\citep{1905.08503} through tropical methods. A survey of some of these approaches can be found in~\citet{maragos2021}. The proper tropical analogue of RKHSs still remained elusive nonetheless.\\
	
	\tb{Main results.} One key property of RKHSs is the versatility entailed by the multiple entry points to their theory, either through a real-valued kernel, a nonlinear feature map, a function space or an integral operator. In this article, we uncover similar links between a tropical kernel, a.k.a.\ a coupling function, its factorization, a tropical function space and monotone operators. Our construction revolves around the famed Fenchel-Moreau conjugations \citep[see e.g.][Chapter 8]{singer1997} with a particular role played by the assumption of symmetry of the kernel. More precisely, as shown by \citet{aronszajn50theory}, the key result in the Hilbertian case is that positive semidefinite kernels coincide with reproducing kernels, i.e.\ 
	\begin{Theorem}[\citet{aronszajn50theory}]\label{thm:aronszajn}
		Given a kernel $k:\X\times \X \rightarrow \R$, the three following properties are equivalent:
		\begin{enumerate}[labelindent=0cm,leftmargin=*,topsep=0.1cm,partopsep=0cm,parsep=0.1cm,itemsep=0.1cm,label=\roman*)]
			\item \label{it:pdk_hilb} $k$ is a \emph{positive semidefinite kernel}, i.e.\ a kernel being both:
			
			\noindent- symmetric: $\forall\, x, y \in X, \; k(x,y)=k(y,x)$, and
			
			\noindent- positive: $\forall\, M\in\N^*,\,\forall \,(a_m,x_m)\in \left(\RR\times X\right)^M,\sum_{n,m=1}^{M} a_n a_mk(x_n,x_m) \geq 0$;	%

			\item \label{it:featMap_hilb} there exists a Hilbert space $(\Hsc,(\cdot,\cdot)_\Hsc)$ and a \emph{feature map} $\Phi:\X\rightarrow \Hsc$ such that 
			
			\noindent- $\forall\, x,y \in X,\, k(x,y)~=~(\Phi(x),\Phi(y))_{\Hsc}$;
			
			\item \label{it:rk_hilb} $k$ is the \emph{reproducing kernel} of the Hilbert space (RKHS) of functions $\Hk:=\overline{\Hsc_{k,0}}$, the completion for the pre-scalar product $( k(\cdot,x),k(\cdot,y))_{k,0}=k(x,y)$ of the space $\Hsc_{k,0}:=\Sp(\{k(\cdot,x)\}_{x\in \X})$, in the sense that
			
			\noindent- $\forall\, x \in X,\; k(\cdot,x)\in\Hk \text{ and }\forall \, f \in \Hsc, \; f(x)=(f,k(\cdot,x))_{\Hsc}$.
		\end{enumerate}
	\end{Theorem}
	Our goal is to obtain the analogue of \Cref{thm:aronszajn} in a max-plus context, and our first main result can be stated informally as follows:
	
	\begin{Theorem}[Tropical analogue of Aronszajn theorem]\label{thm:aronszajn_trop}	Given a kernel $b:\X\times \X \rightarrow \R\cup\{-\infty\}$, the three following properties are equivalent
		\begin{enumerate}[labelindent=0cm,leftmargin=*,topsep=0.1cm,partopsep=0cm,parsep=0.1cm,itemsep=0.1cm,label=\roman*)]
			\item \label{it:tpsd_informal} $b$ is a \emph{tropically positive semidefinite kernel};

			\item \label{it:featMap_informal} there exists a factorization of $b$ by a \emph{feature map} $\psi:\X\rightarrow \pmR^\Z_{\max}$ for some set $\Z$;
			
			\item \label{it:rk_hilb_informal} $b$ is the \emph{sesquilinear reproducing kernel} of a max-plus space of functions $\Rg(B)$, the max-plus completion of $\{\sup_{n\in \{1,\dots,N\}} a_n+b(\cdot, x_n) \,|\, N\in\N^*, a_n\in \mR, x_n\in\X \}$, and $b$ defines a tropical Cauchy-Schwarz inequality over $\pmRX$.
		\end{enumerate}
	\end{Theorem}
	The precise statements are to be found in~\Cref{prop:tpsd_feat_map}, \Cref{prop:largM_CS_cyclic}, the definitions and notation being detailed below. The factorization as in \Cref{thm:aronszajn_trop}-\ref{it:featMap_informal} provides an extension to the infinite dimensional setting of a theorem of \citet{Cartwright2012}. Indeed our analysis tackles general sets $\X$, and we leverage the finite case to foster intuitions, and mostly for counter-examples. Besides, while \citet{aronszajn50theory} stresses the connection of reproducing kernels with inner products as in \Cref{thm:aronszajn}-\ref{it:featMap_hilb}, \citet[][Propositions 8 and 9]{schwartz1964sev} on the contrary approaches kernels through duality pairings. This duality viewpoint will prove important in the following and is summarized in a dictionary to be found in \Cref{tab:hpd_tpsd_concepts}, translating with care the concepts from the Hilbertian to the tropical world.
	
	As a matter of fact, the tropical world can be equipped with two sets of operations, the ``linear'' $(\max,+)$ and the ``sesquilinear'' $(\max,-)$. We shall see, as hinted at by the sesquilinear Legendre-Fenchel transform, that the ``linear'' interpretation of the reproducing property is more restrictive than the ``sesquilinear'' version. Indeed, we show in~\Cref{thm:carac_idempotent} that the ``linear'' reproducing kernels are characterized by von Neumann regularity; in particular, their ranges are images	of linear retractions (idempotent linear maps). It should be stressed that von Neumann regularity is much more restrictive in the tropical setting than in the Hilbertian setting: whereas any closed subspace of a Hilbert space is the range of a linear	retraction (the orthogonal projection over this space), ranges	of tropically linear retractions are rare, and for instance	they do not include spaces of convex functions. This comes as a negative result on the descriptive ability of ``linear'' approaches, studied by e.g.\ \citet{litvinov2011trop}. On the other hand we characterize in \Cref{thm:sym_anti-involution} the ranges of symmetric kernels as the complete submodules with anti-involutions.
	
	Irrespective of this distinction, similarly to Hilbertian kernels, optimization over ranges of tropical kernels enjoys ``representer theorems'', such as \Cref{cor:representer-thm}, ensuring that for some infinite dimensional regression problems, there exist solutions lying in finite dimensional spaces. A prime example of application is convex regression. On the other hand, we prove in \Cref{lem:tpsd_Maup} that choosing a Lagrangian and a set of trajectories defines an idempotent tropical kernel, which is positive definite whenever the Lagrangian is nonnegative, and that this Maupertuis kernel acting on spacetime points generates the value functions. We finally illustrate the representer theorem by an application to inverse optimal control, in which one infers an unknown stopping cost from pointwise measurements of the value function.\\
	
	\paragraph{\bf Related work.} Conjugations such as the Legendre-Fenchel transform are at the core of convex optimization. They also appear in the dual formulation of optimal transport \citep[][Section 4.1]{santambrogio2017}. Indeed costs in optimal transport correspond to the opposite of our kernels since we favour a $(\sup,-)$ convention over the $(\inf,-)$ of optimal transport. The monographies of \citet[Chapter 3.3]{Rachev1998} and \citet[Chapter 2.4]{Villani2003} discuss known results on $c$-concave functions. For instance, the notion of $c$-cyclical monotonicity, further elaborated in \citet[Chapter 5]{Villani2009}, is a key property of optimal transport plans. However optimal transport theory puts little emphasis on the role of the set of $c$-concave functions, which is the ``range'' $\Rg(B)$ of a conjugation $B$ (see \Cref{rmk:opt_transport} for more comments). On the other hand, conjugations have been known since \citet{singer1984conj} to correspond to a tropical kernel \citep[see also][for the extension to general Galois connections]{Akian04setcoverings}. Recent applications of conjugations were investigated for instance in \citet{volle2013,chancelier2021capra}. We here combine the relation between conjugations and kernels with insights from \citet[][Chapter 1]{Pallaschke1997} on monotonicity. The Fenchel transform was itself characterized by \citet{ArtsteinAvidan2009} as the only anti-involution, up to affine change of coordinates, over the set of lower semicontinuous convex functions.
	
	Tropical positive semidefinite matrices have been investigated for discrete sets $\X$ \citep{yu2014tropicalizing,Cartwright2012}. In particular, Yu defined the tropical positive matrices as the image by the nonarchimedean valuation of the set of positive semidefinite matrices over a real closed nonarchimedean field, and showed that they are characterized by the positivity of their $2\times 2$ tropical minors. The tropicalization of the subclass of totally positive matrices was also investigated in~\cite{adiniv}. %
	
	Applying tropical kernels to control problems has two main sources of inspiration. It was recently shown in \citet{aubin2020Riccati,aubin2020hard_control} that, for linear-quadratic optimal control, the Lagrangian cost function and the dynamical system encode the relevant Hilbertian kernel of controlled trajectories. This suggested to incorporate the time component in the variables. Our approach is also inspired by the ``fundamental solutions'' of \citet{dower2015zhang} who consider a tropical kernel on state space, with the time acting as a fixed parameter.\\
	
	It has been often hinted at that max-plus analysis can be interpreted as a limit case of log-sum-exp operations, a property known as Maslov's dequantization \citep{litvinov2005maslov}, akin to taking the limit in Planck's constant in passing from quantum to classical mechanics. Similar interpolating relations between optimal transport and RKHSs have been discussed recently in \citet{feydy2019interpolating}. Our efforts are thus directed toward bridging the gap between the max-plus and the Hilbertian worlds, defining the adequate tropical analogue to RKHSs. This first enquiry opens many theoretical questions, concerning for instance the analogue to the topological characterization of RKHSs or to the minimal factorization of the kernel \citep[see][Theorem 4.21]{steinwart08support}. It also opens up the question of applications of the theory and of the computational advantages of tropical positive semidefiniteness.
	
	The paper is structured as follows. Preliminaries on conjugations are summarized in Section~\ref{sec:preliminaries}. Section~\ref{sec:tpsd_kernel} introduces tropical positive definite kernels and their factorization. In Section~\ref{sec:cyclicity}, these kernels are identified with monotone operators and shown to lead to a tropical Cauchy-Schwarz inequality. In Section~\ref{sec:sesquilinear} we characterize the tropical ranges of symmetric kernels, which we define as tropical reproducing kernel spaces. The ``linear'' approach is investigated in Section~\ref{sec:linear}. We conclude stating a representer theorem in Section~\ref{sec:applications} and applying it to optimization problems, in particular some requiring to find a value function.%

	\section{Tropical functional analysis preliminaries}\label{sec:preliminaries}
	
	\tb{Notations:} Let $\N=\{0,1,\ldots\}$, $\N^* = \{1,2,\ldots\}$, $\R_+$ and $\R_+^*$ denote the set of natural numbers, positive integers, non-negative and positive reals, respectively. We use the shorthand $[M]=\{1,\dots,M\}$. The extended real line is denoted by $\pmR=[-\infty,+\infty]$, we also use the notations $\pR=(-\infty,+\infty]$ and $\mR=[-\infty,+\infty)$. These sets are equipped with the upper (resp.\ lower) addition and subtraction \citep[][p.3]{moreau1970infconvol}, extending the usual operations by indicating through a dot whether $+\infty$ or $-\infty$ is absorbing, e.g.\ $ \infty \dotbplus (-\infty) = -\infty$, $ \infty \dotminus \infty = +\infty$. When $\pmR$ is equipped with the $(\max,+)$ (resp.\ $(\min,+)$) operations we denote it by $\pmR_{\max}$ (resp.\ $\pmR_{\min}$). Given any sets $\X$ and $E$, we denote by $E^\X$ the set of functions $f$ from $\X$ to $E$. If $E$ is equipped with $(\min,+)$ operations we denote generically its elements by $\hat{f}$. A set $\G$ is a complete submodule of $\pmRXmax$ if it is stable under arbitrary sups and addition of constants. When $\X$ is equipped with a topology, a function $f:\X\rightarrow \pmR$ is said to be lower semicontinuous (l.s.c.) if its epigraph is a closed subset of $\X\times\pmR$.\\ %
	
	In all that follows, given a set $\X$, a \emph{kernel} $b$ is a function from $\X\times\X$ to $\pmR$. We shall be interested in at least three specific kernels over $\X=\R^N$ equipped with its Euclidean norm $\|\cdot\|_2$ and inner product $(\cdot,\cdot)_2$ or with a distance $d$:
	\begin{equation}\label{eq:special_kernels}
		\bc(x,y)= (x,y)_2,\quad \bsc(x,y)=-\|x-y\|_2^2,\quad \bL(x,y)=-d(x,y).
	\end{equation}
	These kernels are respectively related to the sets of proper convex l.s.c.\ functions, 1-semiconvex l.s.c.\ functions, and 1-Lipschitz functions \wrT\ the distance $d$ (see p.\pageref{rmk:examples_RgB}). Given a kernel $b$, we also consider the max-plus linear $B$ and sesquilinear $\bB$ operators, defined over $\pmRX$ as
	\begin{align}\label{eq:trop_operators}
		Bf(x)= \sup_{y\in\X} b(x,y) \dotbplus f(y), \quad \bB f(x)= \sup_{y\in\X} b(x,y) \dotbminus f(y), \, \forall\, x\in\X, \, f\in\pmRX,
	\end{align}
	both with the convention that $-\infty$ is absorbing. For instance, for $b=\bc$, $\bB$ is the Fenchel transform which is indeed sesquilinear as $\bB(\min(f,g))=\max(\bB f,\bB g)$ and $\bB(f-\lambda)=\bB f + \lambda$, effectively turning $(\min,-)$ into $(\max,+)$ operations. More formally:
	
	\begin{Definition}\label{def:rmax-sesqui-lin}
		A map $B:\pmRX\rightarrow \pmRX$ is said to be
		\begin{enumerate}[labelindent=0cm,leftmargin=*,topsep=0.1cm,partopsep=0cm,parsep=0.1cm,itemsep=0.1cm,label=\roman*)]
			\item $\pmRmax$-linear if $B(\sup\{f_i\}_{i\in I})=\sup\{Bf_i\}_{i\in I}$ and $B(f \dotplus \lambda)=Bf \dotplus\lambda$ (with $+\infty$ absorbing on both sides), for any finite index set $I$ and $\lambda\in\pmR$; we say in addition that $B$ is continuous if
			$B(\sup\{f_i\}_{i\in I})=\sup\{Bf_i\}_{i\in I}$ holds even for infinite families.
			\item $\pmRmax$-sesquilinear if $B(\inf\{f_i\}_{i\in I})=\sup\{Bf_i\}_{i\in I}$ and $B(f \dotplus \lambda)=Bf\dotbminus\lambda$ (with $+\infty$ absorbing on the l.h.s.\ and  $-\infty$ absorbing on the r.h.s.), for any finite index set $I$ and $\lambda\in\pmR$; 
			we say in addition that $B$ is continuous if
			$B(\inf\{f_i\}_{i\in I})=\sup\{Bf_i\}_{i\in I}$ holds even for infinite families.
		\end{enumerate}
		The range $\Rg(B)$ of $B$ is defined as the set of functions $g\in \pmRXmax$ such that $g=Bf$ for some $f\in \pmRXmax$.
	\end{Definition}
	The continuity notion we introduced needs actually to be defined more generally for our analysis. Recall
	than an ordered set $D$ is directed if for all $d,d'\in D$, there
	exists $d''\in D$ such that $d\leq d''$ and $d'\leq d''$.
	\begin{Definition}\label{def-scott}
		Let $\mathcal{G}$ and $\mathcal{H}$ be complete lattices. A
		map $B:\mathcal{G}\to \mathcal{H}$ will be said to be continuous, if for all directed subsets $D\subset \mathcal{G}$, $B(\sup D) = \sup B(D)$.
	\end{Definition}
	
	This is precisely the notion of continuity with respect to the {\em Scott topology}~\citep{continuous}, which, as noted in~\cite{akiantams,Cohen2004}, is a canonical one in idempotent analysis. In particular, a $\pmRmax$-linear map is an order preserving self-map of $\pmRXmax$, and one can check it is continuous iff the relation $B(\sup\{f_i\}_{i\in I})=\sup\{Bf_i\}_{i\in I}$  holds for all infinite families of functions $\{f_i\}_{i\in I}\subset \pmRX$. Dually, we shall think of a $\pmRmax$-sesquilinear map $B$ as being order preserving, from $\pmRXmax$ equipped with the {\em opposite} to the standard order,  to $\pmRXmax$  equipped with the standard order.
	In this way a $\pmRmax$-sesquilinear $B$ can be seen to be continuous iff the relation $B(\inf\{f_i\}_{i\in I})=\sup\{Bf_i\}_{i\in I}$  holds for all infinite families of functions $\{f_i\}_{i\in I}\subset \pmRX$, in which the infima and suprema refer to the standard order. Observe that no ambiguity will arise from our convention because a (Scott) continuous map is automatically order preserving.
	
	We define the indicator functions $\delta^\bot_{x},\delta^\top_{x}\in \pmRX$ as follows\footnote{Introducing two indicator functions is necessary because of change of signs and of the $\max-\min$ duality.}
	\begin{equation}\label{eq:def_Dirac}
		\delta^\bot_{x}\left(y\right):=\left\{\begin{array}{ll}
			0 & \text { if } y=x, \\
			-\infty & \text { otherwise,} \end{array}\right.
		\delta^\top_{x}\left(y\right):=\left\{\begin{array}{ll}
			0 & \text { if } y=x, \\
			+\infty & \text { otherwise.}
		\end{array}\right.
	\end{equation}
	The $\pmRmax$-sesquilinear and continuous maps, a.k.a.\ \emph{(Fenchel-Moreau) conjugations}, have been characterized as the ones having a kernel as in \eqref{eq:trop_operators}:
	\begin{Proposition}[Theorem 3.1, \citet{singer1984conj}]\label{thm:singer}
		A map $\bB:\pmRX\rightarrow \pmRX$ is $\pmRmax$-sesquilinear and continuous if and only if there exists a kernel $b:\X\times \X \rightarrow \pmR$ such that $\bB f(x)= \sup_{y\in\X} b(x,y) \dotbminus f(y)$. Moreover in this case $b$ is uniquely determined by $\bB$ as $b(\cdot,x)=\bB\delta^\top_{x}$.
	\end{Proposition}
	
	\begin{Remark}\label{rmk:singer_linear} By a change of sign, this result also holds for $\pmRmax$-linear operators, $B f(x)= \sup_{y\in\X} b(x,y)\dotbplus f(y)$.	We refer to \citet[Theorem 2.1]{Akian04setcoverings} for extensions of \Cref{thm:singer} which is the tropical analogue of Riesz representation theorem in Hilbert spaces \citep[see also][]{martinezlegaz1990dualities}. The $\pmRmax$-(sesqui)linear and continuous maps verify that $\Rg(B)$ is a max-plus completion in the sense that it is the smallest complete submodule of $\pmRXmax$ containing $\{\sup_{n\in \{1,\dots,N\}} a_n\dotbplus b(\cdot, x_n) \,|\, N\in\N^*, a_n\in \mR, x_n\in\X \}$
		and %
		\begin{equation}\label{eq:RgB_tpsd}
			\Rg(B)=\{\sup_{x\in \X} a_x\dotbplus b(\cdot, x) \,|\, a_x\in \mR \}.
		\end{equation}
	\end{Remark}

	In the following, we also extensively use the following duality product over $\pmRXmin \times \pmRXmax$, denoting by $\hat{g}$ the elements of $\pmRXmin$,
	\begin{equation}\label{eq:duality_product}
		\langle \hat{g}, f\rangle := \sup_{x\in\X} f(x)\dotbminus \hat{g}(x) \quad \forall (\hat{g},f)\in\pmRXmin \times \pmRXmax.
	\end{equation}
	This duality product allows to define the adjoint of an operator $B:\pmRX\rightarrow \pmRX$.
	
	\begin{Definition}\label{def:selfadjoint}
		If it exists, the adjoint map $\bB'$ of a $\pmRmax$-sesquilinear map $\bB:\pmRX\rightarrow \pmRX$ is defined as the one such that
		\begin{equation*}
			\langle \hat{g}, \bB \hat{f}\rangle = \langle \hat{f}, \bB' \hat{g}\rangle, \quad	\forall (\hat{g},\hat{f})\in\pmRXmin \times \pmRXmin
		\end{equation*}
		If $\bB'=\bB$, then $\bB$ is said to be $\pmRmax$-hermitian. If $\bB$ is continuous with kernel $b(x,y)$, then $\bB'$ exists and corresponds to $b(y,x)$ \citep[Theorem 8.4]{singer1997}.
	\end{Definition}

	\begin{table}[]
		\caption{Corresponding concepts between Hilbertian and tropical kernels}\scriptsize
		\label{tab:hpd_tpsd_concepts}
		\begin{tabular}{@{}c|c|c|c@{}}
			Concept & Hilbertian kernel & Tropical kernel & Reference \\ \hline
			symmetry & $k(x,y)=k(y,x)$  & $b(x,y)=b(y,x)$ & Def.~\ref{def:tpsd_kernel} \\ \hline
			positivity & $\sum_{i,j}a_i a_j k(x_i,x_j)\ge 0$  & $b(x,x)+b(y,y)\ge b(x,y)+b(y,x)$ & Def.~\ref{def:tpsd_kernel} \\ \hline
			feature map & $k(x,y)= (\Phi(x), \Phi(y))_{\Hsc}$ &
			$b(x,y) = \sup_{z\in \Z} \psi(x,z) + \psi(y,z)$ & Prop.~\ref{prop:tpsd_feat_map}  \\ \hline
			\begin{tabular}{c}duality\\ bracket\end{tabular} & $\langle \mu, f \rangle_{\R^{\X,*}\times\R^X}=\int_\X f(y) \d\mu(y)$  & $\langle \hat{g}, f\rangle = \sup_{x\in\X} f(x)-\hat{g}(x)$ & \eqref{eq:duality_product} \\ \hline
			\begin{tabular}{c}kernel\\ operator\end{tabular} & $K(\mu)(x)=\int_\X k(x,y) \d\mu(y)$  & $\bB (\hat{f})(x)= \sup_{y\in\X} b(x,y) \dotbminus \hat{f}(y)$ & Prop.~\ref{thm:singer} \\ \hline
			\begin{tabular}{c}monotone\\operator\end{tabular} & $\langle \mu, K(\mu) \rangle_{\R^{\X,*}\times\R^X}\ge 0$  & $\langle \hat{f}, \bB \hat{f}\rangle \dotplus \langle \hat{g}, \bB \hat{g}\rangle \ge \langle \hat{f}, \bB \hat{g}\rangle \dotbplus \langle \hat{g}, \bB \hat{f}\rangle$ & Prop.~\ref{prop:largM_CS_cyclic} \\[1mm]\hline
			\begin{tabular}{c}function\\ space\end{tabular} & $\Hk=\overline{\Sp(\{k(\cdot,x)\}_{x\in\X})}$  & $\Rg(B)=\{\sup\limits_{x\in \X} [a_x+b(\cdot, x)] \,|\, a_x\in \mR \}$ & Prop.~\ref{thm:singer}+\eqref{eq:RgB_tpsd} \\ \hline
			\begin{tabular}{c}reproducing\\ property\end{tabular}	& $f(x)=(k(\cdot,x),f(\cdot))_{\Hk}$ & $\hat{g}(x)= \langle \bB\hat{g}, \bB \delta^\top_x\rangle=(\bB \bB \hat{g}) (x)$ & Def.~\ref{def:RKMS_rep_prop} %
		\end{tabular}
	\end{table}
	
	\section{Tropical positive semidefinite kernels}\label{sec:tpsd_kernel}
	
	\begin{Definition}\label{def:tpsd_kernel} We say that a kernel $b:\X\times\X\rightarrow \mR$ is a \emph{tropical positive semidefinite} (tpsd) kernel if it is
		\begin{enumerate}[labelindent=0cm,leftmargin=*,topsep=0.1cm,partopsep=0cm,parsep=0.1cm,itemsep=0.1cm,label=\roman*)]
			\item symmetric: $\forall x, y \in \X, \; b(x,y)=b(y,x)$, and
			\item tropically positive: $\forall x, y \in \X, \; b(x,x)+b(y,y)\ge b(x,y)+b(y,x)$.
		\end{enumerate}
	\end{Definition}
	Notice that all the three kernels of \eqref{eq:special_kernels} are tropically positive semidefinite and finite-valued. Moreover \Cref{def:tpsd_kernel} implies that every Hilbertian positive semidefinite kernel is tropically positive semidefinite.\footnote{This is also true for Hilbertian conditionally positive semidefinite (cpsd) kernels, those for which \Cref{thm:aronszajn}-\ref{it:pdk_hilb} only holds for $(a_m)_{m\in[M]}\in \R^M$ such that $\sum_{m=1}^M a_m=0$. Note in the tropical setting we only need to require the property to hold for $M=2$. For $M=1$, tropical positivity is always satisfied since $b(x,x)\ge -\infty$. For $M>2$ we refer to \Cref{prop:largM_CS_cyclic}.} Any square Monge matrix\footnote{A Monge matrix $B\in\pmR^{M\times N}$ is one for which $b_{ij}+b_{mn}\le b_{in}+b_{mj}$, for all $1\le i <m\le M $ and $1\le j <n\le N$. They correspond to submodular functions over the discrete set $\X \times \X$, we refer to \citet{Burkard1996} for a review of their properties and to \citet{Weiss2016} for their application to the assignment problem.} corresponds to a tpsd kernel, since the latter relaxes the Monge requirements on all the 2-by-2 minors to only the principal minors. Similarly, the logarithm of the absolute value of any positive semidefinite kernel is a tpsd kernel.\footnote{This result is related to the relation between Hilbertian cpsd and indefinitely divisible psd kernels \citep[see e.g.][]{Berg1984}. For positive semidefinite kernels, the Gram matrices $\begin{pmatrix}
			k(x,x) & k(x,y)\\
			k(x,y) & k(y,y)
		\end{pmatrix}$ are positive semidefinite. Thus, when considering their determinant, we obtain that $k(x,x)k(y,y)\ge k(x,y)k(y,x)$. Taking the logarithm (since we allow tpsd kernels to have the $-\infty$ value), we obtain the tropical positivity of $\log(|k|)$.} An immediate consequence of \Cref{def:tpsd_kernel} is that the set of tropically positive semidefinite kernels is stable for the sum, when adding a constant in $\mR$, for pointwise limits, and when restricting the kernel to subsets of $\X$. Other examples of kernels include the opposite of powers of distances over any metric space, such as the Wasserstein distances in optimal transport, or the kernel $b(X,Y)=-\|\operatorname{Spec} \log(XY^{-1})\|^2_2$ over positive semidefinite matrices acting on $\R^N$, with $\operatorname{Spec}$ denoting the eigenspectrum.%

	We chose to consider $\mR$-valued kernels even though most of our kernels of interest are $\R$-valued. Nevertheless, the Dirac mass $\delta^\bot_x$ is $\mR$-valued and plays the role of a neutral element. Note that $\delta^\bot_x$ is the logarithm of the Hilbertian psd ``0-1'' kernel defined as $k(x,x)=1$, $k(x,y)=0$ if $x\neq y$. It is possible to extend our analysis to $\pmR$-valued kernels, although it would be cumbersome to keep under check two infinite values in our computations, as stressed in \citet[Remark 8.25b]{singer1997}. We first show that all tpsd kernels are translations of symmetric diagonal-vanishing and nonpositive-valued kernels (which can be interpreted as costs in a game context, we will return to this in \Cref{sec:least-action}). This property may remind of \citet[Lemma 2.1, Proposition 3.2]{Berg1984} stated for negative definite kernels. %
	\begin{Lemma}\label{lem:tpsd_negative-valued}
		A kernel $b: \X \times \X \to \mR$ is tpsd if and only if there exists a function $\phi:\X\rightarrow \mR$ and a symmetric kernel $b_0: \X \times \X \to \mR$, with $b_0(x,x)= 0$ and $b_0(x,y)\le 0$ for all $x,y\in\X$, such that
		\begin{equation}\label{eq:tpsd_negative-valued}
			b(x,y) =  \phi(x)+b_0(x,y)+\phi(y).
		\end{equation}
		Moreover, given $b(\cdot,\cdot)$, we have that $\phi(x)=b(x,x)/2$ and $\Rg(B)=\phi+\Rg(B_0)$.
	\end{Lemma}
	\begin{proof}
		If $b$ is tpsd, then we set $\phi(x)=b(x,x)/2$. If $(\phi(x),\phi(y))\in\R^2$, we pose $b_0(x,y):=b(x,y)-\phi(x)-\phi(y)$, and set it to $0$ otherwise. If $\phi(x)=-\infty$ for some $x\in\X$, then, by the tropical positivity of $b$, $b(x,y)=-\infty$ for all $x\in\X$, whence $b_0$ satisfies \eqref{eq:tpsd_negative-valued}. Again the positivity of $b$ yields that $b_0$ is symmetric, diagonal-vanishing and nonpositive-valued. Conversely, if $b$ is defined by \eqref{eq:tpsd_negative-valued}, it is trivially tropically positive since $b_0$ vanishes on the diagonal and is nonpositive-valued. Let $f\in\Rg(B)$, as per \eqref{eq:RgB_tpsd}, we can find $(a_x)_{x\in\X}\in \mR^\X$ such that, for any $y\in\X$,
		\begin{equation*}
			f(y)=\sup_{x\in \X} a_x+b(y, x)=\phi(y)+\sup_{x\in \X} a_x+\phi(x)+b_0(y, x)
		\end{equation*}
		so $f\in\phi+\Rg(B_0)$, the converse is shown similarly.
	\end{proof}

	We now draw the connection with the characterization of Hilbertian positive semidefinite kernels $k$ through feature maps, as in \Cref{thm:aronszajn}-\ref{it:featMap_hilb}. The result below shows the analogy with such a characterization for the tropical inner product over $\mR^\Z$ defined by $(f,g)_{\sup}=\sup_{z\in\Z} f(z)\dotbplus g(z)$.
	\begin{Proposition}\label{prop:tpsd_feat_map}
		Let $b: \X \times \X \to \mR$ be a kernel.
		The following properties are equivalent
		\begin{enumerate}[labelindent=0cm,leftmargin=*,topsep=0.1cm,partopsep=0cm,parsep=0.1cm,itemsep=0.1cm,label=\roman*)]
			\item\label{it_tpsd} $b$ is tpsd; %
			\item\label{it_featMap} there exists a set $\Z$ and a function $\psi: \X\times \Z \to \mR$ such that
			\begin{equation}\label{eq:feature_map_carac}
				b(x,y) = \sup_{z\in \Z} \psi(x,z) + \psi(y,z).
			\end{equation}
		\end{enumerate}
	\end{Proposition}
	\begin{proof}
		\ref{it_tpsd}$\Rightarrow$\ref{it_featMap}.
		Take $\Z=\X\times \X$, and consider the function $\psi$ such that, for all $x,y\in \X\times \X$, $\psi(x,(x,y))=b(x,x)/2$ and $\psi(x,(y,x))=b(x,y)\dotbminus b(y,y)/2$, with $-\infty$ absorbing. We set to $-\infty$ all the other values of $\psi(x,(u,v))$, for which $x\not \in \{u,v\}$, whence $\psi$ takes its values in $\mR$. Then, by definition of $\psi$, since $-\infty$ is absorbing, the only $z$ for which the values $\psi(x,z)$ and $\psi(y,z)$ can be finite are $z\in\{(x,y),(y,x)\}$. As $b$ is symmetric, we obtain that
		\begin{align*}
			\sup_{z\in \Z} \psi(x,z) \dotbplus \psi(y,z)\\
			&\hspace{-3cm}= \max \left(\psi(x,(x,y)) \dotbplus \psi(y,(x,y)), \psi(x,(y,x)) \dotbplus \psi(y,(y,x))\right)\\
			&\hspace{-3cm}= \max\left(\frac{b(x,x)}{2}+b(y,x)\dotbminus\frac{b(x,x)}{2}, b(x,y)\dotbminus\frac{b(y,y)}{2}+\frac{b(y,y)}{2}\right)=b(y,x), 
		\end{align*}
		the equality holding for $b(x,x)=-\infty$ or $b(y,y)=-\infty$, since $b(y,x)<\infty$ by assumption.
		
		\ref{it_featMap}$\Rightarrow$\ref{it_tpsd}. The kernel $b$ as defined by \eqref{eq:feature_map_carac} is symmetric and we assumed it takes its values in $\mR$. Then, for all $z\in \Z$,
		\begin{equation*}
			2(\psi(x,z) + \psi(y,z))\le 2 \sup_{z\in \Z} \psi(x,z)+ 2\sup_{z\in \Z}\psi(y,z)= b(x,x)+b(y,y),
		\end{equation*} 
		taking the supremum over $z\in \Z$ allows us to conclude. 
	\end{proof}
	
	\noindent \tb{Examples of factorizations $\psi(x,z)$ with $\Z=\X$}:	Note that if \eqref{eq:feature_map_carac} holds for $\Z=\X$ and $\psi=b$, then $-b$ satisfies a triangular inequality. For instance, for the kernels defined in \eqref{eq:special_kernels}, it is true for $\bL$, but it does not hold for $\bc$ or $\bsc$, for which the feature map cannot be as simple as $\psi=b$.
	\begin{enumerate}[label=\roman*),labelindent=0em,leftmargin=2em,topsep=0cm,partopsep=0cm,parsep=0cm,itemsep=2mm]
		\item For $\X=\R^N$ and $\bc(x,y)=(x,y)_2$, we can choose $\psi(x,z)=\frac{1}{2}\|x\|^2_2-\|x-z\|_2^2$.
		\item For $\X=\R^N$ and $\bsc(x,y)=-\|x-y\|_2^2$, we can choose $\psi(x,z)=-2\|x-z\|_2^2$.
		\item For $(\X,d)$ a metric space, $b(x,y)=-d(x,y)^p$ with $p\in(0,1]$, we can choose $\psi=b$ as a consequence of the subadditivity of $t\in\R_+\mapsto t^p$.\\
	\end{enumerate}
	
	However we cannot always take $\Z=\X$, even for finite sets. Indeed, the issue of the factorization of $b$ has been studied in~\citet{Cartwright2012} when $\X$ is a finite set, of cardinality $|\X|=n$. Then, the minimal cardinality
	of $Z$ such that there exists a factorization \eqref{eq:feature_map_carac} (summarized by the notation $b=\psi\psi'$) is called the tropical symmetric Barvinok rank. It is shown in \citet[Theorem~4]{Cartwright2012}, that, for every tpsd $b$, there exists a factorization \eqref{eq:feature_map_carac} with $|\Z|\leq \max(n,\lfloor n^2/4\rfloor)$, and that this bound is tight. Hence there exists kernels $b$ over finite sets $\X$ for which all $\Z$ verify $|\Z|\ge\lfloor n^2/4\rfloor)>n$. For instance for $n=5$, defining the kernel through its Gram matrix, we set
	\begin{displaymath}
		B=\begin{pmatrix}
			0 & \text{-}1 & 0 & 0 & 0\\
			\text{-}1 & 0 & 0 & 0 & 0\\
			0 & 0 & 0 & \text{-}1 & \text{-}1\\
			0 & 0 & \text{-}1 & 0 & \text{-}1\\
			0 & 0 & \text{-}1 & \text{-}1 & 0
		\end{pmatrix}.
	\end{displaymath}
	Such $B$ corresponds to the complete bipartite graph $K_{\floor{\frac{n}{2}}, \ceil{\frac{n}{2}}}$, for which the smallest cardinality of $\Z$ as in \eqref{eq:feature_map_carac} is the clique cover number, equal to $6>n=5$ \cite[][Remark 3.3]{Cartwright2012}.

	Nevertheless the feature map $\Phi_b(x):= b(\cdot,x)$ remains of interest even if it does not correspond to a factorization \eqref{eq:feature_map_carac}. As a matter of fact, the notion of convexity \wrT\ the family of functions $\Phi_b(\X)=\{b(\cdot,x)\}_{x\in\X}$ has been extensively studied in \citet{Pallaschke1997} and \citet[Chapter 9]{singer1997}. In particular the tropical SDP inequality corresponds to the monotonicity of $\Phi_b$ in the sense that, for any $x,y\in\X$ \citep[][p.5, (1.1.8)]{Pallaschke1997}
	\begin{equation*}
		\Phi_b(x)(x)+\Phi_b(y)(y)-\Phi_b(x)(y)-\Phi_b(y)(x)=b(x,x)+b(y,y)-b(x,y)-b(y,x)\ge 0.
	\end{equation*}
	Monotonicity of $\Phi_b$ makes it a candidate to be the subgradient of a $\Phi_b$-convex function. In Banach spaces, maximal monotone operators are obtained as subgradients of convex functions. In our tropical setting, while we actually do not have the desired maximality, we show in \Cref{sec:cyclicity} that tropical positivity entails monotonicity and that we have an analogue of the Cauchy-Schwarz inequality. This will require to switch our view to the operator $\bB$ corresponding to $b$.%

	\section{Monotonicity of operators induced by tpsd kernels}\label{sec:cyclicity}
	
	The following result explains why, unlike for the Hilbertian kernels of \Cref{thm:aronszajn}-\ref{it:pdk_hilb}, considering only pairs, rather than sequences of points, is sufficient in the tropical setting. \Cref{prop:largM_CS_cyclic} furthermore relates the tropical positive semidefiniteness of the kernel over elements of $\X$ to the (cyclic) monotonicity of the operator $\bB$ over functions of $\pmRX$.

	\begin{Proposition}\label{prop:largM_CS_cyclic}
		Given a kernel $b:\X\times\X\rightarrow \mR$, set $\bB$ as in \eqref{eq:trop_operators}. Then the following statements are equivalent:
		\begin{enumerate}[label=\roman*),labelindent=0cm,leftmargin=*,topsep=0.1cm,partopsep=0cm,parsep=0.1cm,itemsep=0.1cm]
			\item \label{it_tpsd2} the kernel $b$ is tpsd, i.e.\ symmetric and positive: $$\forall x, y \in \X, \; b(x,x)+b(y,y)\ge b(x,y)+b(y,x);$$
			\item \label{it_tpsdM} for all $M\in \N^*$, $(x_m)_{m\in[M]}\in \X^M$ and permutations $\sigma:[M]\rightarrow [M]$, the Gram matrices $\b G:=[b(x_n,x_m)]_{n,m\in[M]}$ are symmetric and satisfy%
			\begin{equation}\label{eq:tpsd_matrices}
				\sum_{m=1}^{M} b(x_m,x_m) \ge \sum_{m=1}^{M} b(x_m,x_{\sigma(m)});
			\end{equation}
			\item \label{it_monot2} the operator $\bB$ is $\pmRmax$-hermitian and monotone for the duality pairing in the sense that
			\begin{equation}\label{eq:monot_operator}
				\forall\,\hat{f}, \hat{g}\in\pmRX,\,  \langle \hat{f}, \bB \hat{f}\rangle \dotplus \langle \hat{g}, \bB \hat{g}\rangle \ge \langle \hat{f}, \bB \hat{g}\rangle \dotbplus \langle \hat{g}, \bB \hat{f}\rangle;
			\end{equation}
			\item \label{it_monot2_max} the operator $\bB$ is $\pmRmax$-hermitian and 
			\begin{equation}\label{eq:monot_operator_max}
				\forall\,\hat{f}, \hat{g}\in\pmRX,\,  \max\left(\langle \hat{f}, \bB \hat{f}\rangle, \langle \hat{g}, \bB \hat{g}\rangle\right) \ge \langle \hat{f}, \bB \hat{g}\rangle;
			\end{equation}
			with, for any $x$, $b(x,x)=-\infty$ implying that for all $y\in\X$, $b(x,y)=-\infty$;
			\item \label{it_monotM} the operator $\bB$ is $\pmRmax$-hermitian and cyclic monotone for the duality pairing in the sense that
			\begin{gather}
				\forall\,M\in \N^*,\,\forall\,(\hat{f}_m)_{m\in[M]}\in\pmR^{\X,M} \text{ with the convention } f_{M+1}=f_1,\nonumber\\
				\sum_{m=1}^{M} \langle \hat{f}_{m}, \bB \hat{f}_{m}\rangle \ge \sum_{m=1}^{M} \langle \hat{f}_{m}, \bB \hat{f}_{m+1}\rangle;\label{eq:cycMonot_operator}
			\end{gather}
			\item \label{it_monotM_max}  the operator $\bB$ is $\pmRmax$-hermitian and 
			\begin{gather}
				\forall\,M\in \N^*,\,\forall\,(\hat{f}_m)_{m\in[M]}\in\pmR^{\X,M} \text{ with the convention } f_{M+1}=f_1,\nonumber\\ \max_{m=1,\dots M} \langle \hat{f}_{m}, \bB \hat{f}_{m}\rangle \ge \max_{m=1,\dots M} \langle \hat{f}_{m}, \bB \hat{f}_{m+1}\rangle,\label{eq:cycMonot_operator_max}
			\end{gather}
			with, for any $x$, $b(x,x)=-\infty$ implying that for all $y\in\X$, $b(x,y)=-\infty$.
		\end{enumerate}
	\end{Proposition}
	\begin{Remark}[Relation with optimal transport]\label{rmk:opt_transport} The notion of cyclic monotonicity is a cornerstone of optimal transport \citep[Chapter 5]{Villani2009}, expressing the fact that there is no gain in permuting the assignments of an optimal transport plan (ibid., Theorem 5.9). In other words, the support of the optimal plan for continuous costs $c$ is concentrated on a $c$-cyclically monotone subset of $X\times X$. Correspondingly, \eqref{eq:tpsd_matrices} characterizes the diagonal $\{(x,x)\}_{x\in\X}$ as one such set whenever $b$ is tpsd. Similarly, \Cref{lem:tpsd_negative-valued} states that tpsd kernels are translations of nonpositive kernels vanishing on the diagonal. It is customary in optimal transport to assume nonnegative costs but \citet[Remark 4.2]{Villani2009} underlines that this is mostly to have a lower-bound on $c$. Besides, symmetry is not paramount and rarely assumed, except in specific cases such as \citet[Section 2.5]{di2017optimal}. For a cost $c$ such that $(-c)$ is tpsd, the optimal transport cost between two measures is, after translation, symmetric, nonnegative and vanishes on the diagonal, so only the triangle inequality is missing to make it a (pseudo)-distance.%
	\end{Remark}	
	
	\begin{proof}\ref{it_tpsdM}$\Rightarrow$\ref{it_tpsd2}. This corresponds to the case $M=2$ and the cycle $(\sigma(1),\sigma(2))=(2,1)$.
		
		\ref{it_tpsd2}$\Rightarrow$\ref{it_tpsdM}. Since every permutation can be decomposed as a product of cycles, we just have to prove the property for the latter. For $M=1$, we have by definition that $b(x,x)\ge -\infty$ for all $x\in\X$. Fix $M\in\N^*$ with $M\ge 2$. Reindexing $(x_m)_{m\in[M]}$, we can assume w.l.o.g.\ that $\sigma(m)=m+1$ for $m\in[M-1]$ and $\sigma(M)=1$. Hence, since the kernel is tpsd, with the convention that $x_{M+1}=x_{1}$
		\begin{align*}
			\sum_{m=1}^{M} b(x_m,x_m)&=\frac{1}{2}\sum_{m=1}^{M} [b(x_m,x_m)+b(x_{m+1},x_{m+1})]\\ &\ge \sum_{m=1}^{M} b(x_m,x_{m+1})=\sum_{m=1}^{M} b(x_m,x_{\sigma(m)}).
		\end{align*}
		
		\ref{it_tpsd2}$\Leftrightarrow$\ref{it_monot2}. Notice that for all $x,y\in\X$, $\langle \delta^\top_x, \bB \delta^\top_y\rangle=b(x,y)$. Hence \ref{it_tpsd2} is just a specialization of \ref{it_monot2}. To prove \ref{it_tpsd2}$\Rightarrow$\ref{it_monot2}, note for every $\hat{f},\hat{g}\in \pmRXmin$, according to \eqref{eq:duality_product} and since $b$ is tpsd, we have that
		\begin{align*}
			\langle \hat{f}, \bB \hat{f}\rangle \dotplus \langle \hat{g}, \bB \hat{g}\rangle&=\sup_{x,u\in\X} b(x,u)\dotbminus\hat{f}(x)\dotbminus\hat{f}(u) \dotplus \sup_{y,v\in\X} b(y,v)\dotbminus\hat{g}(y)\dotbminus\hat{g}(v)\\
			&\ge b(x,x)+ b(y,y) \dotbminus 2 \hat{f}(x)\dotbminus 2\hat{g}(y) \text{ setting $u=x, v=y$}\\
			&\ge 2 b(x,y) \dotbminus 2 \hat{f}(x)\dotbminus 2\hat{g}(y).
		\end{align*}
		Taking the supremum over $x,y\in\X$ yields the result.	
		
		\ref{it_monot2}$\Leftrightarrow$\ref{it_monotM}. This equivalence follows from the same computations as for the proof of \ref{it_tpsd2}$\Leftrightarrow$\ref{it_tpsdM} Cyclic monotonicity is precisely the $\sigma$-cycle considered in the proof above.
		
		\ref{it_monot2}$\Rightarrow$\ref{it_monot2_max}. This stems from the fact that, applying \eqref{eq:monot_operator} and the hermitianity of $\bB$,
		\begin{equation*}
			2\max\left(\langle \hat{f}, \bB \hat{f}\rangle, \langle \hat{g}, \bB \hat{g}\rangle\right)\ge \langle \hat{f}, \bB \hat{f}\rangle \dotplus \langle \hat{g}, \bB \hat{g}\rangle \ge \langle \hat{f}, \bB \hat{g}\rangle \dotbplus \langle \hat{g}, \bB \hat{f}\rangle = 2 \langle\hat{f}, \bB \hat{g}\rangle.
		\end{equation*}
		\ref{it_monot2_max}$\Rightarrow$\ref{it_tpsd2}. Since $\bB$ is $\pmRmax$-hermitian, $b$ is symmetric. Fix $x,y\in\X$. If $b(x,x)+ b(y,y)\neq-\infty$, then we can define the functions $\hat{f}(\cdot)=\delta^\top_x(\cdot)-\frac{b(x,x)}{2}$ and $\hat{g}(\cdot)=\delta^\top_y(\cdot)-\frac{b(y,y)}{2}$ and apply \eqref{eq:monot_operator_max}, whence
		\begin{equation*}
			0=\max\left(\langle \hat{f}, \bB \hat{f}\rangle, \langle \hat{g}, \bB \hat{g}\rangle\right)\ge \langle\hat{f}, \bB \hat{g}\rangle=b(x,y)-\frac{b(x,x)}{2}-\frac{b(y,y)}{2}.
		\end{equation*}
		If $b(x,x)+ b(y,y)=-\infty$, then, by \ref{it_monot2_max}, $ b(x,y)=-\infty$. In both cases, we have shown that $b$ is tropically positive.
		
		\ref{it_monot2_max}$\Leftrightarrow$\ref{it_monotM_max} We just have to consider $M=2$ for \ref{it_monotM_max}$\Rightarrow$\ref{it_monot2_max}, while \ref{it_monot2_max}$\Rightarrow$\ref{it_monotM_max} is a consequence of the transitivity of the maximum.
	\end{proof}
	
	The two formulas \eqref{eq:monot_operator}-\eqref{eq:monot_operator_max} can be seen as two different ways of defining a tropical analogue of the monotonicity of a nonlinear operator over a Hilbert space, informally written as ``$\langle \hat{f}-\hat{g}, K \hat{f}- K \hat{g}\rangle \ge 0$''. We can moreover interpret \eqref{eq:monot_operator_max} as a tropical Cauchy-Schwarz inequality, informally written as ``$\langle \hat{f}, K \hat{f}\rangle^{\frac{1}{2}} \langle \hat{g}, K \hat{g}\rangle^{\frac{1}{2}} \ge \langle \hat{f}, K \hat{g}\rangle$'' From the notion of monotonicity \eqref{eq:monot_operator}, we can define a form of tropical ``quadratic'' discrepancy $d_B:\pmRXmin\times \pmRXmin \rightarrow \pmR$:
	\begin{align}
		d_B(\hat{f},\hat{g}):=\frac{1}{2}\left[\langle \hat{f}, \bB \hat{f}\rangle \dotplus \langle \hat{g}, \bB \hat{g}\rangle \dotminus \langle \hat{f}, \bB \hat{g}\rangle \dotminus \langle \hat{g}, \bB \hat{f}\rangle\right]. \label{eq:dB} 
	\end{align}
	\tb{Examples of $d_B$:} For  a metric space $(\X,d)$ and $b(x,y)=-d(x,y)$, we shall see below (p.\pageref{rmk:examples_RgB}) that $\Rg(B)$ is the set of 1-Lipschitz functions (including the constant functions in $\pmR$) and that $\bB \hat{f}=-\hat{f}$ for $\hat{f}\in\Rg(B)$. Hence, for any $\hat{f},\hat{g}\in\Rg(B)$, we have that $\langle \hat{g}, \bB \hat{f}\rangle= - \inf_\X (\hat{f}+\hat{g})$ and $2d_B(\hat{f},\hat{g})= \inf_\X (\hat{f}+\hat{g})-\inf_\X \hat{f}-\inf_\X\hat{g}$. So $d_B(\hat{f},\hat{g})=0$ indicates that $\hat{f}$ and $\hat{g}$ have a common infimum. A similar computation for $b(x,y)=\delta^\bot_x(y)$ gives the same $d_B$, but defined over $\pmRXmin$. Nevertheless, for non-compact $\X$, even with smooth kernels, $d_B$ can be infinite. Indeed for $b(x,y)=(x,y)_2$, $ \langle b(\cdot,w), \bB b(\cdot,w)\rangle= \sup_{x,y\in\R^d} (x,y)_2 \dotbminus (w,x+y)_2=\infty$. Note that \Cref{prop:largM_CS_cyclic} could be rewritten in a $\pmRXmax$-linear rather than $\pmRXmax$-sesquilinear setting, replacing the duality product by a sup-inner product and changing signs, $\langle \hat{f}, \bB \hat{g}\rangle=\sup_{x,y\in\X} b(x,u)\dotbminus\hat{f}(x)\dotbminus\hat{g}(y)=:(-\hat{f},B(-\hat{g}))_{\sup}$.

	\section{A sesquilinear theory: defining Reproducing Kernel Moreau Spaces based on symmetry}\label{sec:sesquilinear}
	
	We shall be interested in spaces that are ranges of tropical operators. We start with a characterization of such ranges when the kernels are symmetric.
	
	\begin{Theorem}\label{thm:sym_anti-involution}
		Let $\G$ be a complete submodule of $\pmRXmax$. Then the following statements are equivalent:
		\begin{enumerate}[label=\roman*),labelindent=0cm,leftmargin=*,topsep=0.1cm,partopsep=0cm,parsep=0.1cm,itemsep=0.1cm]
			\item \label{it_rgB} there exists a symmetrical kernel $b:\X\times \X \rightarrow \pmR$ such that $\G=\Rg(B)$;
			\item \label{it_antiinvol} there exists a $\pmRmax$-sesquilinear map $\bF:\G\rightarrow \G$ such that $\bF \bF=\Id_{\G}$, i.e.\ $\bF$ is an anti-involution over $\G$.
		\end{enumerate}
		If these properties hold, then $\bF$ can be taken as the restriction of $\bB$ to $\Rg(B)$.
	\end{Theorem}
	Here we only focus on the existence and do not discuss the uniqueness of the anti-involution. For classical convexity, \citet{ArtsteinAvidan2009} characterized the Fenchel conjugation as the only anti-involution, up to an affine change of coordinates,  over the set of l.s.c.\ convex functions $\G$. They furthermore interpret the condition of \Cref{thm:sym_anti-involution}-\ref{it_antiinvol} as defining a ``duality'' over $\G$ (ibid., Definition 11). %
	
	\begin{proof} \ref{it_rgB}$\Rightarrow$\ref{it_antiinvol}. We use the classical property that $\bB \bB \bB = \bB$, see e.g.\  \citep[][p.3]{Akian04setcoverings}, so $\bB$ is an anti-involution over $\Rg(B)$.

		\ref{it_antiinvol}$\Rightarrow$\ref{it_rgB}. We first show that $\bF$ is necessarily continuous in the sense of~\Cref{def-scott} and that it can be extended to a $\pmRmax$-hermitian and continuous map $\tilde{F}:\pmRX\rightarrow \pmRX$.

		Given an arbitrary index set $\Asc$, fix a family $(g_\alpha)_{\alpha\in \Asc}\in\G^{\Asc}$. We define its infimum relatively to $\G$ as
		\[ \inf\nolimits^\G_\alpha g_\alpha:=\max\{g\in \G \, |\, \forall\alpha\in \Asc, \, g\le g_{\alpha}\}  ,\]%
		noting that the latter set does admit a greatest element, as $\G$ is a complete submodule of $\pmRXmax$. In particular, if the family $(g_\alpha)_{\alpha\in\Asc}$ consists
		of a single element $g$, then,
		\[ \inf\nolimits^\G g
		= \max\{h \in \G\,|\, h\leq g
		\}  .
		\]
		The operation of relative supremum with respect to $\G$ is defined in a dual manner, however, since $\G$ is stable by arbitrary suprema, it merely coincides with the ordinary supremum of functions. As $\bF$ is $\pmRmax$-sesquilinear, it satisfies that
		\begin{equation*}
			\forall f,g\in\G \text{ s.t.\ } f\le g, \bF(f)=\bF(\min(f,g))=\max(\bF(f),\bF(g))\ge \bF(g),
		\end{equation*}
		in other words, $\bF$ is antitone. Hence, as $\inf\nolimits^\G_{\alpha\in \Asc} g_\alpha \le g_\beta$ for any $\beta\in \Asc$,
		\begin{displaymath}
			\bF(\inf\nolimits^\G_{\alpha\in \Asc} g_\alpha)\ge \sup_{\alpha\in \Asc}\bF( g_\alpha) \ge \bF( g_\beta).
		\end{displaymath}
		Since $\G$ is a complete submodule of $\pmRXmax$, $\sup_{\alpha\in \Asc}\bF( g_\alpha)\in\G$, so composing by $\bF$ and using that $\bF$ is antitone, we derive that
		\begin{equation*}
			\bF\bF(\inf\nolimits^\G_{\alpha\in \Asc} g_\alpha)\le \bF(\sup_{\alpha\in \Asc}\bF( g_\alpha)) \le \bF\bF( g_\beta).
		\end{equation*}
		As $\bF$ is an anti-involution, taking the infimum in $\G$ over $\beta$ on the r.h.s., yields
		\begin{equation*}
			\inf\nolimits^\G_{\alpha\in \Asc} g_\alpha = \bF\bF(\inf_{\alpha\in \Asc} g_\alpha) \le \bF(\sup_{\alpha\in \Asc}\bF( g_\alpha)) \le \inf\nolimits^\G_{\beta\in \Asc} g_\beta
			\enspace.
		\end{equation*}
		So $\inf\nolimits^\G_{\alpha\in \Asc} g_\alpha = \bF\left(\sup_{\alpha\in \Asc}\bF( g_\alpha)\right)$, and, composing again by $\bF$, we deduce that $\bF$ is continuous on its domain, meaning that $\bF$ sends relative infima with respect to $\G$ to relative suprema (which coincide with ordinary suprema).
		
		Define now $\tilde{F}:\pmRX\rightarrow \pmRX$ by, for any $f\in\pmRX$, $\tilde{F}(f):=\bF(\inf\nolimits^\G f) \in \G$. The function $\tilde{F}$ is automatically an extension of $\bF$. We now show that $\tilde{F}$ is $\pmRmax$-hermitian. Fix $f\in\pmRX$ and $\lambda\in\pmR$, then, as $\bF$ is $\pmRmax$-sesquilinear, $\tilde{F}(f\dotplus \lambda)=\bF(\lambda\dotplus\inf\nolimits^\G f)=\tilde{F}(f)\dotbminus\lambda$. Let $(f_\alpha)_{\alpha\in \Asc}\in(\pmRX)^{\Asc}$. By definition, $\inf\nolimits^\G(\inf_{\alpha\in \Asc} f_\alpha)=\inf\nolimits^\G_\alpha (\inf\nolimits^\G f_\alpha)$. As $\bF$ is continuous, 
		\begin{equation*}
			\tilde{F}(\inf_{\alpha\in \Asc} f_\alpha) = \bF(\inf\nolimits^\G_\alpha (\inf\nolimits^\G f_\alpha))=\sup_\alpha \bF(\inf\nolimits^\G f_\alpha)= \sup_\alpha\tilde{F}(f_\alpha). 
		\end{equation*}
		Hence $\tilde{F}$ is $\pmRmax$-hermitian and continuous. We can now apply \Cref{thm:singer} to derive a kernel $b$ associated with $\tilde{F}$. By construction, $\Rg(\tilde{F})\subset\G$ and since $\bF$ is an anti-involution, $\G=\Rg(\bF)\subset \Rg(\tilde{F})$, so $\G=\Rg(\tilde{F})$. As $\bF$ is an anti-involution, 
		\begin{displaymath}
			\tilde{F}\tilde{F}(f)=\bF\bF(\inf\nolimits^\G f)=\inf\nolimits^\G f\le f.
		\end{displaymath}
		We have thus recovered one of the characterizations of dual Galois connections over $\pmRXmax$, as given in \citep[][p.3, Eq.(2a)]{Akian04setcoverings}, so $\tilde{F}$ is equal to its dual connection. \citet[Theorem 2.1, Proposition 2.3]{Akian04setcoverings} then allows to conclude that, since the dual connection is simply the transpose $\tilde{F}^\top$ \citep[see also][Theorem 8.4]{singer1997}, the kernel $b$ is symmetric.
	\end{proof}

	\begin{Remark}
		\Cref{thm:sym_anti-involution} should be compared with
		\citet[Th.~23]{DS04} and \citet[Th.~42]{Cohen2004}, which entail
		that the row and column spaces of a tropical matrix
		are anti-isomorphic lattices. By specializing this
		result, we deduce that the row space of a symmetric
		matrix is anti-isomorphic to itself. \Cref{thm:sym_anti-involution}
		refines this result (and also extends it to the infinite dimensional setting), showing that a version of the latter anti-isomorphism property
		characterizes the ranges of symmetric operators.
	\end{Remark}
	We now leverage the characterization, $\bB \bB \bB \hat{f} = \bB \hat{f}$, expressed for symmetric kernels in \Cref{thm:sym_anti-involution}, to define tropical sesquilinear reproducing kernel spaces.	As shown in \eqref{eq:RgB_tpsd}, using spaces $\Rg(B)$ to define RKMSs provides a direct analogy with the fact that a RKHS $\Hk$ is the completion for its norm $\|\cdot\|_k$ of $\Sp(\{k(\cdot,x)\}_{x\in\X})$.
	
	\begin{Definition}\label{def:RKMS_rep_prop} We call \emph{reproducing kernel Moreau spaces} (RKMS) the complete submodules $\Rg(B)$ of $\pmRXmax$  where $\bB$ is a $\pmRmax$-sesquilinear continuous and hermitian operator associated with the symmetric kernel $b$. For all $\hat{g}\in \Rg(B)$ and $x\in\X$, we say that they satisfy a \emph{sesquilinear reproducing property}
		\begin{align}
			\hat{g}(x)&= (\bB \bB \hat{g}) (x)= \langle \bB\hat{g}, \bB \delta^\top_x\rangle= \sup_{z\in\X} b(z,x) \dotbminus [\sup_{y\in\X} b(z,y) \dotbminus \hat{g}(y)]. \label{eq:repr-prop_b}
		\end{align}
	\end{Definition}  
	The sesquilinear reproducing property as defined in \eqref{eq:repr-prop_b} is not an empty statement. It characterizes the elements of $\G=\Rg(B)$ through an immediate lemma, proved again using the identity $\bB \bB \bB= \bB$.
	
	\begin{Lemma}\citep[Corollary 8.5]{singer1997}\label{lem:RKMS_rep_prop} Let $\bB$ be a $\pmRmax$-sesquilinear and continuous operator. Then for any $g\in\pmRX$, $\hat{g}=\bB \bB \hat{g}$ holds if and only if $g\in \Rg(B)$.
	\end{Lemma}		
	For $b=\bc$, $\bB$ is the Fenchel conjugate operator, whence \eqref{eq:repr-prop_b} is equivalent to Fenchel's theorem stating that convex l.s.c.\ functions are the only fixed points of the Fenchel biconjugate. Note that the difficult part in Fenchel's theorem is to identify $\Rg(B)$, i.e.\ proving that all convex l.s.c.\ functions are outputs of the Fenchel transform.
	
	We choose to interpret $\hat{g}=\bB \bB \hat{g}$ as a reproducing property, however $\bB \bB$ is not a $\pmRXmax$-(sesqui)linear operator. One may thus wonder whether there is a more direct interpretation. Indeed we can relate \eqref{eq:repr-prop_b} to a $\pmRXmin$-linear operator, but it unfortunately cannot characterize $\Rg(B)$. This will be further emphasized in \Cref{sec:linear}.
	
	\begin{Lemma}\label{lem:Funk_RKMS} For all $\hat{g}\in \Rg(B)$ and $x\in\X$, \eqref{eq:repr-prop_b} is equivalent to
		\begin{align}
			\hat{g}(x)&=\inf_{y\in\X} \hat{g}(y) \dotplus \sup_{z\in\X}[b(z,x)- b(z,y)]=: \inf_{y\in\X} \hat{g}(y)\dotplus  c(x,y)=:(\Cop \hat{g}) (x) \label{eq:repr-prop_c}
		\end{align}
		where $c(x,y):=\sup_{z\in\X}[b(z,x) \dotbminus b(z,y)]$ is the Funk distance between $b(\cdot,x)$ and $b(\cdot,y)$, and $\Cop:\pmRXmin\rightarrow \pmRXmin$  the related $\pmRmin$-linear operator.	
	\end{Lemma}
	\begin{proof}
		The term $c(x,y)$ appears by permuting the `$\sup\inf$` into `$\inf\sup$` and using \eqref{eq:repr-prop_b} since
		\begin{align*}
			\hat{g}(x)&=\sup_{z\in\X} b(z,x) \dotbminus [\sup_{y\in\X} b(z,y) \dotbminus \hat{g}(y)]=\sup_{z\in\X}\inf_{y\in\X} b(z,x) \dotbminus b(z,y) \dotplus \hat{g}(y)\\
			&\le \inf_{y\in\X} \hat{g}(y)\dotplus \sup_{z\in\X} [b(z,x) \dotbminus b(z,y)]= \inf_{y\in\X} \hat{g}(y) \dotplus c(x,y)\stackrel{y=x}{\le} \hat{g}(x).
		\end{align*}
	\end{proof}
	Unfortunately $c$ does not correspond to a unique $b$, nor to a unique $\Rg(B)$, as shown in the examples below:
	
	\noindent \tb{Examples of tpsd $b(x,y)$, $c(x,y)$ and $\Rg(B)$}:
	\begin{enumerate}[labelindent=0em,leftmargin=1.5em,topsep=0cm,partopsep=0cm,parsep=0cm,itemsep=2mm,label=\roman*)]\label{rmk:examples_RgB}
		\item For $\X=\R^N$, $b(x,y)=(x,y)_2$ gives $c(x,y)=\delta^\top_x(y)$ whereas $\Rg(B)$ is the set of proper convex l.s.c.\ functions adding the constant functions $\pm\infty$ \citep[][Theorem 3.7]{singer1997}.
		\item For $\X=\R^N$, $b(x,y)=-\|x-y\|^2$ gives $c(x,y)=\delta^\top_x(y)$ whereas $\Rg(B)$ is the set of proper 1-semiconvex l.s.c.\ functions adding the constant functions $\pm\infty$ \citep[][Theorem 3.16]{singer1997}.
		\item For any $\X$ and $\alpha\ge 0$, $b(x,y)=\left\{\begin{array}{ll}
			0 & \text { if } y=x, \\
			-\alpha & \text { otherwise,}	\end{array}\right.$ gives $c(x,y)=-b(x,y)$ whereas $\Rg(B)$ is the set of functions $f$ which difference $f(x)-f(y)$ is smaller than $\alpha$. For $\alpha=+\infty$, $b(x,y)=\delta^\bot_x(y)$, $\Rg(B)$ corresponds to the whole $\pmRX$ \citep[][Remark 3.2]{singer1997}.
		\item For $(\X,d)$ a metric space and $p\in(0,1]$, $b(x,y)=-d(x,y)^p$ gives $c(x,y)=d(x,y)^p$ whereas $\Rg(B)$ is the set of $(1,p)$-Hölder continuous functions \wrT\ the distance $d$ (i.e.\ $|f(x)-f(y)|\le 1\cdot d(x,y)^p$), when adding the constant functions $\pm\infty$ \citep[][Theorem 3.14]{singer1997}.\\
	\end{enumerate}
	
	\noindent \tb{Informal analogy between RKMS and RKHS}: We can draw an informal analogy between \eqref{eq:repr-prop_c} and the reproducing property of functions $f$ of a RKHS $\Hk$ defined through a differential operator $D$ (and its adjoint $D^*$) over a bounded open set $\X$, i.e.\
	\begin{align}
		f(x)&= (f(\cdot),k(\cdot,x))_{\Hk}\overset{(*)}{=}\int_{y\in\X} f(y) D_y^*D_y k(x,y)\d y \label{eq:repr-prop_k}
	\end{align}
	where $(*)$ holds for Green kernels with null integral on the boundaries, as per \citet[Section 1.7]{saitoh16theory}. Equation \eqref{eq:repr-prop_k} underlines the fact that the integral was formally replaced with an inf, the product with a sum and the differential operator $D^*D$ by the sup of a difference, all this in a classical tropical fashion. Moreover for the RKHSs defined through differential operators, \eqref{eq:repr-prop_k} expresses the fact the functions $f\in\Hk$ are fixed points for a kernel integral operator, but whose kernel is not simply $k(x,y)$.\\
	
	\noindent \tb{Reproducing but not positive}: By \Cref{thm:sym_anti-involution}, any kernel $b(x,y)$ that is symmetric, even non-tpsd, gives $\bB \bB \bB = \bB$, so the reproducing property \eqref{eq:repr-prop_b} can be defined for spaces and kernels that do have not have a notion of positivity. This is true also for Hilbertian kernels, since reproducing properties can be defined for indefinite kernels and their related Krein spaces \citep[see e.g.][and references within]{ong2004learning}. Even the symmetry requirement could be relaxed as was done for tropical operators by \citet{Akian04setcoverings,singer1997} and for RKHSs by \citet{mary2005theory}.\\

	\noindent \tb{Relation between $\bB$ and $\Cop$}:
	The tropical positive operators $\bB$ are $\pmRXmax$-sesquilinear over $\pmRX$ whereas $\Cop$ is $\pmRXmin$-linear. Based on the non-symmetric $c(x,y)$ one can define the set of $c$-Lipschitz functions.
	\begin{equation*}
		\Lipb(\X,\pmR):=\{g\in\R^{\pm,\X}\, \mid\, g(x)\le g(y)\dotplus c(x,y)\}.
	\end{equation*}
	For finite-valued kernels $b$, since $c(x,x)=0$, $\Lipb(\X,\pmR)$ corresponds to the fixed points ($g=\Cop g$) of $\Cop$. By \eqref{eq:repr-prop_b},
	\begin{equation*}
		\Rg\nolimits^{\max}(B)\subseteq \Lipb(\X,\pmR) \subseteq \Rg\nolimits^{\min}(\Cop):=\{\inf_{x\in \X} a_x+c(\cdot, x) \,|\, a_x\in \pR \}
	\end{equation*}
	with equality when $c(x,y)=-b(x,y)$ (i.e.\ $-b$ is idempotent). Notice that $\Rg\nolimits^{\max}(B)$ is sup-stable, $\Lipb(\X,\pmR)$ is sup-stable and inf-stable while $\Rg\nolimits^{\min}(\Cop)$ is inf-stable.%

	\section{A linear theory: only idempotent operators and spaces of Lipschitz functions}\label{sec:linear}
	
	We have seen in \Cref{sec:sesquilinear} that one could introduce some tropically linear, rather than sesquilinear, operators. This choice can be related to the identification of kernel integral operators in the Hilbertian context with $\pmRXmax$-linear operators in the tropical setting as was done by \citet[Section 7,][]{litvinov2011trop}. We prove below that discarding sesquilinearity imposes considerable restrictions on the types of spaces that one can consider. These are either ranges of idempotent kernels or lattices of Lipschitz functions. This result further emphasizes why we need to consider $\pmRmax$-sesquilinear operators to properly define the max-plus analogue of RKHSs. Recall for instance that convex l.s.c.\ functions are defined based on the $\pmRmax$-sesquilinear Fenchel transform, not a linear version of it.\\
	
	For RKHSs $(\Hk,(\cdot,\cdot)_{\Hk})$, the reproducing property writes as follows
	\begin{equation}\label{eq:rep_prop_hilb}
		\forall\, f\in\Hk, \,\forall\, x\in\X,\; f(x)=\langle \delta^{lin}_{x}, f(\cdot) \rangle_{\R^{\X,*}\times\R^\X}=(k(\cdot,x),f(\cdot))_{\Hk},
	\end{equation}
	where $k(\cdot,x)\in \Hk$ and the space of finite measures $\R^{\X,*}$ is the dual of $\R^\X$ for the pointwise convergence. The latter is generated by the linear Dirac masses $\delta^{lin}_{x}:f\in\R^\X\mapsto f(x)$.  %
	
	Consider a set $\G\subset\pR^\X$. %
	A $\pmRXmax$-linear analogy to \eqref{eq:rep_prop_hilb} would be the existence of a kernel $c:\X\times \X\rightarrow \pmR$ satisfying
	\begin{equation}\label{eq:rep_prop_trop}
		\forall\, g\in\G, \,\forall\, x\in\X,\; g(x)=\sup_{y\in\X} \delta^\bot_{x}(y) \dotbplus g(y) = \sup_{y\in\X}c(x,y) \dotbplus g(y).
	\end{equation}
	Although it is always possible to consider $c(x,y)=\delta^\bot_{x}(y)$ in \eqref{eq:rep_prop_trop}, since $\G\subset\pR^\X$, $\delta^\bot_{x}(\cdot)\notin\G$. Nonetheless among the possible $c$, there is one that is maximal:
	\begin{Proposition}\label{prop:cG_kernel}
		For any set $\G\subset\pR^\X$, there is a kernel $c_\G$ that is maximal over the kernels $c$ satisfying \eqref{eq:rep_prop_trop}, for the partial order inherited from $\pmR$. Moreover $c_\G:\X\times \X\rightarrow \pmR$ is given by
		\begin{equation}\label{eq:maximal_kernel}
			c_\G(x,y):=\inf_{g\in\G} g(x) \dotminus g(y).
		\end{equation}
		If, for some $x\in\X$, there exists a function $g\in\G$ such that $g(x)<\infty$, i.e.\ $\G$ is \emph{proper} at $x$, then $c_\G(x,x)=0$ (otherwise  $c_\G(x,x)=+\infty$). %
		Moreover $c_\G$ is idempotent in the $(\max,+)$ sense
		\begin{equation}\label{eq:gammaG_idem}
			c_\G(x,y) = \sup_{z\in \Z} c_\G(x,z) \dotbplus c_\G(z,y).
		\end{equation}
	\end{Proposition}
	\begin{proof} We first show that $c_\G$ satisfies \eqref{eq:rep_prop_trop}. Fix any $x\in\X$ and $g\in\G$, and define
		\begin{equation}\label{eq:maximal_kernel_rep_prop}
			\tilde{g}(x)= \sup_{y\in\X} g(y)\dotbplus c_\G(x,y)=\sup_{y\in\X} g(y)\dotbplus\inf_{h\in\G} h(x)  \dotminus h(y)\stackrel{h=g}{\le} g(x).
		\end{equation}%
		If $h(x)=+\infty$ for all $h\in\G$, since $+\infty$ is absorbing, $c_\G(x,y)=+\infty$ for all $y\in\X$. So $\tilde{g}(x)\ge g(x)+c_\G(x,x)=\infty+\infty=\infty=g(x)$. Assume now that $\{g\in\G, g(x)<\infty\}\neq \emptyset$. Since $+\infty$ is absorbing and $\G\subset\pR^\X$, $c_\G(x,x)=\inf_{g\in\G, g(x)<\infty} g(x) \dotminus g(x)=0$. Hence, by \eqref{eq:maximal_kernel_rep_prop}, $g(x)\ge\tilde{g}(x)\ge g(x)\dotplus c_\G(x,x)=g(x)$, whence $g=\tilde{g}$. By definition of $\tilde{g}$ in \eqref{eq:maximal_kernel_rep_prop}, we deduce that $c_\G$ satisfies \eqref{eq:rep_prop_trop}. Consider now any kernel $c$ satisfying \eqref{eq:rep_prop_trop}, whence, for all $g\in\G$ and all $x,y\in\X$,
		\begin{equation*}
			g(x)\ge  g(y)\dotbplus c(x,y).
		\end{equation*}
		Taking the infimum over $g\in\G$, by definition of $c_\G$ \eqref{eq:maximal_kernel}, $c_\G$ is indeed larger than any $c$. 
		
		Let $x\in\X$. If $h(x)=+\infty$ for all $h\in\G$, then $c_\G(x,y)=+\infty$ for all $y\in\X$, so \eqref{eq:gammaG_idem} holds. Otherwise, we can fix $f\in\G$ such that $f(x)<\infty$. Let $z\in\X$. If $f(z)<\infty$, then
		\begin{equation}\label{eq:proof_cG_idem}
			c_\G(x,z) \dotbplus c_\G(z,y)=\inf_{g\in\G} g(x) \dotminus g(z) \dotbplus \inf_{h\in\G} h(z) \dotminus h(y) \stackrel{g=h=f}{\le} f(x)\dotminus f(y).
		\end{equation}
		Assume now that $f(z)=\infty$. As $f(x)<+\infty$, $c_\G(x,z)=-\infty$ which is absorbing for the l.h.s.\ of \eqref{eq:proof_cG_idem}. Hence $c_\G(x,z) + c_\G(z,y)=-\infty\le  f(x)\dotminus f(y)$. We are thus allowed to take the infimum of the r.h.s.\ of \eqref{eq:proof_cG_idem} over all $f\in\G$, which yields $c_\G(x,z) + c_\G(z,y)\le c_\G(x,y)$. Since, for $\G$ proper at $x$, we also have 
		\begin{equation*}
			\sup_{z\in \Z} c_\G(x,z) \dotbplus c_\G(z,y) \stackrel{z=x}{\ge}c_\G(x,y),
		\end{equation*}
		we deduce that \eqref{eq:gammaG_idem} holds, which concludes the proof.
	\end{proof}
	
	Alike $\Cop$ in \Cref{sec:sesquilinear}, the main limitation of $c_\G$ is that it cannot characterize uniquely each $\G\subset\pR^\X$. Indeed we show in \Cref{prop:cG_asClosure} below that $\Rg(C_\G):=\{\sup_{x\in\R^N} [c_\G( \cdot, x)+a_x] \,|\, a_x\in \mR, x\in\X \}$ is the ``closure'' in a max-min sense of $\G$ and is hence inf and sup-stable.
	\begin{Proposition}\label{prop:cG_asClosure}
		Let $\G\subset\pR^\X$ be a set of $\pmR$-proper functions, and $c_\G$ as in \eqref{eq:maximal_kernel}. Consider $C_\G:\pmRXmax\rightarrow \pmRXmax$ as in \eqref{eq:trop_operators}. Then $C_\G\ge \Id_{\pmRX}$, in the sense that $C_\G(f)(x)\ge f(x)$ for all $f\in \pmRX$, and $\Rg(C_\G)$ is the smallest, in the inclusion sense, inf-stable complete submodule of $\pmRXmax$ containing $\G$. It corresponds to the set of $(\text{-}c'_\G)$-Lipschitz functions, $\Rg(C_\G)=\{f\in\pmRX\,|\, \forall\, x,y\in\X,\, f(x)\le f(y)\dotminus c_\G(y,x)\}$ where $c'_\G$ is the transpose of $c_\G$.
	\end{Proposition}
	\begin{proof}
		Let $\F$ be an inf-stable complete submodule of $\pmRXmax$ containing $\G$. Hence $c_\G(\cdot,y)=\inf_{g\in\G} g(\cdot) \dotminus g(y) \in \F$ since $g(\cdot) \dotminus g(y)\in \F$. The set $\F$ is also sup-stable, so $\Rg(C_\G)\subset \F$. Since $c_\G(x,x)\in\{0,\infty\}$ by \Cref{prop:cG_kernel}, $C_\G(f)(x)= \sup_{y\in\X} c_\G(x,y) \dotbplus f(y)\stackrel{y=x}{\ge} f(x)$ for all $f\in \pmRX$.
		
		We now prove that $\Rg(C_\G)$ is an inf-stable complete submodule of $\pmRXmax$. The set $\Rg(C_\G)$ is by definition sup-stable. Let us show that it is also inf-stable. Take any family $(f_\alpha)_{\alpha\in \Asc}$ with $f_\alpha\in\pmRX$ and $\Asc$ an arbitrary index set. Using successively that $C_\G=C_\G C_\G$ as shown in \Cref{prop:cG_kernel} (as $\G$ is a set of proper functions), then a permutation of inf and sup, and that $C_\G\ge \Id_{\pmRXmax}$, we derive that
		\begin{align*}
			\inf_{\alpha\in \Asc} C_\G f_\alpha = \inf_{\alpha\in \Asc} C_\G C_\G f_\alpha&= \inf_{\alpha\in \Asc} \sup_{y\in\X} C_\G f_\alpha(y)\dotbplus c_\G(\cdot,y)\\
			&\hspace{-1cm}\ge  \sup_{y\in\X} \inf_{\alpha\in \Asc} C_\G f_\alpha(y)\dotbplus c_\G(\cdot,y)=C_\G(\inf_{\alpha\in \Asc} C_\G f_\alpha) \ge \inf_{\alpha\in \Asc} C_\G f_\alpha,
		\end{align*}
		so $\inf_{\alpha\in \Asc} C_\G f_\alpha=C_\G(\inf_{\alpha\in \Asc} C_\G f_\alpha)\in \Rg(C_\G)$ which is thus inf-stable. Finally, as $c_\G$ satisfies \eqref{eq:rep_prop_trop} by \Cref{prop:cG_kernel}, for any $g\in\G$, $g=C_\G g\in \Rg(C_\G)$, so $\G\subset \Rg(C_\G)$.%
		
		Finally we show that $\Rg(C_\G)$ is the set of $(\text{-}c'_\G)$-Lipschitz functions. Let $f$ be a $(\text{-}c'_\G)$-Lipschitz function. By definition, for any $x,y\in\X$, $f(x)\dotbplus c_\G(y,x) \le f(y)$, with equality when $x=y$, so  $f(y)=\sup_x c_\G(y,x) \dotbplus f(x) $, and $f\in \Rg(C_\G)$. Take now $f= C_\G g \in \Rg(C_\G)$, thus we obtain that%
		\begin{align*}
			C_\G g(x)-C_\G g(y)&=\sup_z c_\G(x,z) \dotbplus g(z) - \sup_w c_\G(y,w) \dotbplus g(w)\\ &\hspace{-2cm}\stackrel{w=z}{\le} \sup_z c_\G(x,z) - c_\G(y,z)\\ &\hspace{-2cm} \stackrel{\eqref{eq:gammaG_idem}}{\le} \sup_z c_\G(x,z) - c_\G(y,x) - c_\G(x,z)=- c_\G(y,x),
		\end{align*}		
		which concludes the proof.

	\end{proof}\label{prop:mcG_kernel}

	What if we had considered a weighted inner product in \eqref{eq:rep_prop_trop}? Then we could indeed obtain a formula akin to a reproducing property but only in a very limited context of ranges of idempotent operators. Notice also that this approach does not require the symmetry of the kernels involved.
	\begin{Theorem}\label{thm:carac_idempotent} Let $\G$ be a complete submodule of $\pmRXmax$. The following properties are equivalent:
		\begin{enumerate}[labelindent=0cm,leftmargin=*,topsep=0.1cm,partopsep=0cm,parsep=0.1cm,itemsep=0.1cm,label=\roman*)]
			\item\label{it_WeightInner} there exists two kernels $a,b:\X\times\X\rightarrow \pmR$ such that, when considering the weighted sup-inner product $(f,g)_{\sup, a}=\sup_{y,z\in\X} f(y) + g(z) + a(y,z)$, then, for all $x\in\X$ and $f\in\G$,
			\begin{equation}\label{eq:weighted_rep_prop}
				f(x)=(f,b(x,\cdot))_{\sup, a} \text{ and } b(\cdot,x)\in\G;
			\end{equation}	
			\item\label{it_RgIdem} $\G=\Rg(P)$ for some linear continuous idempotent operator $P: \pmRXmax\rightarrow \pmRXmax$, i.e.\ $P=P^2$;
			\item\label{it_RgVN} $\G=\Rg(B)$ for some von Neumann regular operator $B: \pmRXmax \rightarrow \pmRXmax$, i.e.\ a linear and continuous operator such that there exists a linear continuous $A: \pmRXmax\rightarrow \pmRXmax$ satisfying $B=BAB$.
		\end{enumerate}
	\end{Theorem}
	\begin{proof}
		\ref{it_WeightInner}$\Rightarrow$\ref{it_RgVN}. Fix $a,b$ as in \ref{it_WeightInner}, then, for all $f\in\G$,
		\begin{equation*}
			f(\cdot)=\sup_{y,z\in\X} b(\cdot,y) + f(z) + a(y,z),
		\end{equation*}	
		so $f$ is an arbitrary supremum of elements of $\Rg(B)$, hence $f\in\Rg(B)$. Conversely, since $b(\cdot,x)\in\G$ for all $x\in\X$, $\Rg(B)\subset\G$, so $\Rg(B)=\G$. Hence \eqref{eq:weighted_rep_prop} is equivalent to say that $f=BA f$ for all $f\in \Rg(B)$, so $B=BAB$ and $B$ is von Neumann regular. 
		
		\ref{it_RgVN}$\Rightarrow$\ref{it_WeightInner}. Since $\G=\Rg(B)$, and $B=BAB$, for all $x\in\X$ and $f\in\G$, \eqref{eq:weighted_rep_prop} is satisfied when applying \Cref{thm:singer} to the $\pmRXmax$-linear and continuous operators $A,B$ to derive the corresponding kernels $a,b$ (see \Cref{rmk:singer_linear}).
		
		\ref{it_RgIdem}$\Rightarrow$\ref{it_RgVN}. Setting $B=P$ and $A=\Id_{\pmRX}$ yields the result. In other words, every $(\max,+)$ idempotent operator is von Neumann regular.
		
		\ref{it_RgVN}$\Rightarrow$\ref{it_RgIdem}. Since $BA=BABA$, setting $P=BA$ proves that $P$ idempotent. Since $B=BAB=PB$, we have that $\Rg(B)\subset \Rg(P)$ and, as $P=BA$, $\Rg(P)\subset \Rg(B)$. Hence $\Rg(P)=\Rg(B)=\G$.
	\end{proof}
	
	\begin{figure}[!ht]
		\begin{center}
			
			\def\coord#1#2#3{{-sqrt(3)/2*(#1-#2)} ,{ -(1/2)*#1 - (1/2)*#2 + #3}}
			
			\begin{tikzpicture}
				
				\coordinate (v1) at (\coord{1}{0}{-1});
				\coordinate (v2) at (\coord{0}{0}{0});
				\coordinate (v3) at (\coord{-1}{0}{1});
				\coordinate (v4) at (\coord{0}{-1}{0});
				\coordinate (v5) at (\coord{0}{-0.5}{0});
				\coordinate (v6) at (\coord{1}{0}{0});
				\coordinate (v7) at (\coord{0}{0}{1});
				
				\coordinate (p1) at (\coord{1}{0}{0});
				\coordinate (p2) at (\coord{0}{0}{1});
				\coordinate (p3) at (\coord{0}{-1}{0});

				\filldraw[gray,draw=black,opacity=0.8,very thick] (v1) -- (p1);
				\filldraw[gray,draw=black,opacity=0.8,very thick] (v3) -- (p2);
				
				\filldraw[gray,draw=black,opacity=0.8,very thick] (p1) -- (p3) -- (p2) -- (v2) -- cycle;

				\draw[dashed,->] (\coord{0}{0}{0}) -- (\coord{3}{0}{0}) node[above] {$x_1$};
				\draw[dashed,->] (\coord{0}{0}{0}) -- (\coord{0}{2}{0}) node[above] {$x_2$};
				\draw[dashed,->] (\coord{0}{0}{0}) -- (\coord{0}{0}{2}) node[above, right] {$x_3$};
				\filldraw[black] (\coord{0}{0}{0}) circle (1.5pt) node[below,right] {0};
				
				\filldraw (v1) circle (0.3ex) node[below] {$a$};
				\filldraw (v3) circle (0.3ex) node[below] {$a'$};
				\filldraw (v2) circle (0.3ex) node[below] {$b$};
				\filldraw (v4) circle (0.3ex) node[anchor=south east] {$b'$};
				\filldraw (v6) circle (0.3ex) node[anchor=south east] {$c$};
				\filldraw (v7) circle (0.3ex) node[anchor=south east] {$c'$};
				
			\end{tikzpicture}
		\end{center}
		\caption{A representation of the range of a tropically PSD matrix that is $\pmRmax$-sesquilinear reproducing (\Cref{def:RKMS_rep_prop}), but not $\pmRmax$-linear reproducing (\Cref{thm:carac_idempotent}), since it is not the range of an idempotent matrix.}\label{fig-flag}
	\end{figure}
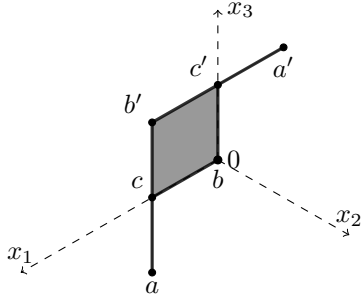
	
	The main drawback of the $\pmRXmax$-linear formula \Cref{thm:carac_idempotent}-\eqref{eq:weighted_rep_prop} is that it does not cover cases of interest which the $\pmRXmax$-sesquilinear \Cref{def:RKMS_rep_prop}-\eqref{eq:repr-prop_b} can tackle. Take for instance $\X=\{-1,0,1\}$, and $b(x,y)=(x,y)_2$, then $\Rg(B)$ is not the range of an idempotent $\pmRXmax$-linear operator. Indeed, the latter are exactly tropical eigenspaces which
	have been characterized in the setting of tropical spectral theory~\cite[Th.~3.100]{bcoq}. The only full-dimensional tropical eigenspaces are indeed
	{\em alcoved polyhedra}, i.e.\ sets defined
	by collections of inequalities $x_i-x_j\geq M_{ij}$ for some matrix
	$M\cup\{-\infty\}$  \citep{AlcovedPolytopes}. Here, the range of $B$, whose
	cross section by the hyperplane orthogonal to the main diagonal
	is shown on~\Cref{fig-flag}, is the union of such
	a three dimensional alcoved polyhedron (the two dimensional cell
	of the cross section), and of two cells
	of lower dimension $2$ (edges of the cross section).
	The points $a,b,a'$, which correspond to the columns of the matrix $B$,
	are the tropical generators of the range of $B$.
	The action of the involution $\bar{F}$ (see~\Cref{thm:sym_anti-involution})
	is illustrated on the figure, it acts as a central symmetry on the parallelogram $(c,b,c',b')$,
	whereas it sends the segment $[c,a]$ to $[c',a']$.
	
	\section{Application to optimization problems}\label{sec:applications}
	\subsection{A representer theorem for general tropical kernels}
	
	One important perk of Hilbertian kernel methods stems from the fact that some infinite-dimensional optimization problems can be solved through equivalent finite-dimensional problems. This behavior is expressed through ``representer theorems'', which ensure that the solutions of an optimization problem live in a finite-dimensional subspace of the RKHS and consequently enjoy a finite representation. These ``representer theorems'' can be informally summarized as: ``a finite number of evaluations implies a finite number of coefficients'' or ``all the information is contained in the samples''. Representer theorems can be found as early as in \citet[Theorems 3.1 and 5.1]{kimeldorf1971tcheb} for quadratic norm-penalties and exact or (quadratic) approximate interpolation \citep[see][for an extension to more general objectives]{scholkopf2001representer}. We show below that tropical kernels enjoy a similar property.
	
	Given two sets $\X$ and $\X'$, a kernel $b:\X\times\X'\rightarrow \mR$, and a subset $\hat{\X}\subset \X$, define
	\begin{multline}\label{eq:def_diffRgB}
		\Rg\nolimits_{\partial\text{-}\hat{\X}}(B):=\Big\{f\in\Rg(B)\,|\, \forall\, x\in\hat{\X},\, \exists\, p_x\in\X', \, f(x)=b(x,p_x)\dotbminus\bB f(p_x) \\
		=\sup_{p\in\X'} b(x,p)\dotbminus\bB f(p)\Big\}.
	\end{multline}
	This set can be understood as the subset of functions of $\Rg(B)$ for which there exists a $B$-subdifferential at every point of $\hat{\X}$.\footnote{We refer to \citet[Section 2]{MartinezLegaz1995subdifferentials} and \citet[Section 3]{Akian04setcoverings} for more comments on this notion and on its relation with the continuity and coercivity of $b$.} For instance, for $\bc(x,y)= (x,y)_2$ and $\X=\X'=\R^N$, it is well-known that $\Rg_{\partial\text{-}\X}(B)$ contains the continuous convex functions, but that is strictly smaller than $\Rg(B)$ since convex l.s.c.\ functions may have an empty subdifferential at points that do not lie in the relative interior of their domain. %
	For functions belonging to $\Rg\nolimits_{\partial\text{-}\hat{\X}}(B)$, \Cref{prop:finite-solution}, below, states that, given a certain number of their values, there always exists another interpolating function represented by the same number of coefficients.\footnote{This may be thought of as an analogue of Carath\'eodory theorem, since, as soon as they are satisfiable by an element $f\in \Rg\nolimits_{\partial\text{-}\hat{\X}}(B)$,  the $\cI$ constraints $y_m=f(x_m)$ can be also satisfied by an element $f$ written as a supremum of at most $\cI$ generators $b(\cdot,p_m)$.  We refer the reader to~\cite{DS04,GK,BSS} for background on the discrete tropical analogue of Carath\'eodory theorem.}
	This entails straightforwardly a representer theorem, \Cref{cor:representer-thm}.

	\begin{Proposition}[Tropical interpolation]\label{prop:finite-solution}
		Let $\cI$ be a nonempty index set, given $(x_m,y_m)_{m\in\cI}\in (\X\times \R)^{\cI}$, setting $\hat{\X}=\{x_m\}_{m\in\cI}$, the three following statements are equivalent:
		\begin{enumerate}[label=\roman*),labelindent=0cm,leftmargin=*,topsep=0.1cm,partopsep=0cm,parsep=0.1cm,itemsep=0.1cm]
			\item \label{it_interp} there exists $f\in \Rg\nolimits_{\partial\text{-}\hat{\X}}(B)$ such that $y_m=f(x_m)$ for all $m\in\cI$;
			\item \label{it_finiteinterp} there exists $(p_m)_{m\in\cI}\in(\X')^{\cI}$ such that $y_m=f^0(x_m)$ for all $m\in\cI$, for
			\[ 
			f^0(\cdot):=\bB(\inf_{m\in\cI}[\delta^\top_{p_m}(\cdot)\dotplus b(x_m,p_m)-y_m])= \max_{m\in\cI} b(\cdot,p_m)\dotbminus b(x_m,p_m) + y_m;
			\]
			\item \label{it_ineqinterp}  there exists $(p_m)_{m\in\cI}\in(\X')^{\cI}$ such that $y_n-y_m\ge b(x_n,p_m)\dotbminus b(x_m,p_m)$ for all $n,m\in\cI$.
		\end{enumerate}
	\end{Proposition}
	\begin{proof}
		\ref{it_interp}$\Rightarrow$\ref{it_finiteinterp}. Set $p_m$ for each $x_m$ as in \eqref{eq:def_diffRgB}, then, for all $n\in\cI$,
		\begin{align*}
			y_n=b(x_n,p_n)\dotbminus\bB f(p_n)&\le \max_{m\in\cI} b(x_n,p_m)\dotbminus\underbrace{(b(x_m,p_m)-y_m)}_{=\bB f(p_m)}=f^0(x_n)\\ &\le\sup_{p\in\X'} b(x_n,p)\dotbminus\bB f(p)=\bB \bB f(x_n)=f(x_n)=y_n,
		\end{align*}
		so $f^0(x_n)=y_n$. 
		
		\ref{it_finiteinterp}$\Rightarrow$\ref{it_interp}. We directly have that $f^0\in \Rg(B)$ and that $p_m$ is a subdifferential at each $x_m$, whence $f^0\in \Rg\nolimits_{\partial\text{-}\hat{\X}}(B)$.
		
		\ref{it_finiteinterp}$\Leftrightarrow$\ref{it_ineqinterp}. This follows from the definition of $f^0$.
	\end{proof}
	
	\begin{Corollary}[Representer theorem]\label{cor:representer-thm}
		Given points $(x_m)_{m\in\cI}\in \X^{\cI}$ and a function $\cL:\pmR^{\cI}\rightarrow \pmR$, fix $\hat{\X}=\{x_m\}_{m\in\cI}$. Then, if the problem
		\begin{equation}\label{eq:opt-infdim}
			\min_{f\in\Rg(B)}\cL((f(x_m))_{m\in\cI})
		\end{equation}
		has a solution $\bar f\in \Rg\nolimits_{\partial\text{-}\hat{\X}}(B)$ with finite values $(f(x_m))_{m\in\cI}\in \R^{\cI}$, it also has a solution $f^0$ as in \Cref{prop:finite-solution}-\ref{it_finiteinterp} which can be obtained solving
		\begin{gather}\label{eq:opt-finDim}
			\min_{(p_m,y_m)_{m\in\cI}\in (\X'\times \R)^M}\cL((y_m))_{m\in\cI})\\ \text{ s.t.\ $y_n-y_m\ge b(x_n,p_m)\dotbminus b(x_m,p_m),\, \forall \, n,m\in\cI$.} \nonumber
		\end{gather}
		Conversely, if \eqref{eq:opt-finDim} has a solution, then it is also a solution in $\Rg\nolimits_{\partial\text{-}\hat{\X}}(B)$ of \eqref{eq:opt-infdim}.
	\end{Corollary}
	When $b$ is the standard scalar product and $\cI$ is finite, each $p_m$ can be interpreted as a subgradient at $x_m$, and \Cref{cor:representer-thm} recovers a well-known property in convex regression \citep[e.g.][Section 6.5.5]{boyd_vandenberghe_2004}, where \eqref{eq:opt-finDim} is then a convex problem provided $\cL$ is convex. We have thus shown that this result also holds for very general kernels and uncountable set $\cI$, not even assuming symmetry or tropical positivity of $b$. Consequently \Cref{prop:finite-solution} should be related to interpolation theorems such as \citet[Theorem 4]{Taylor2016}. For instance, for $\mu$-strongly convex functions with $L$-bounded gradient, this corresponds to $\X=\R^d$, $\X'=\{p\in\R^d\,\|p\|_2\le L\}$ and $b(x,p)=(x,p)_2+\mu \|x\|_2^2$. Besides, if one considers a family of kernels $(b_\alpha)_{\alpha\in \Asc}$ parametrized by some index set $\Asc$ then one can also optimize on the kernel, addressing the problem $\min_{\alpha\in\Asc, f\in\Rg(B_\alpha)}\cL(\alpha,(f(x_m))_{m\in\cI})$. Examples of such families for $\alpha\in \R_+$ are the kernels $b_\alpha(x,y)=-\alpha\|x-y\|_2^2$, corresponding to $\alpha$-semiconvex functions, or $b_\alpha(x,y)=-\alpha d(x,y)$, corresponding to $\alpha$-Lipschitz functions. This adds a further difficulty, as, for instance, when optimizing over $\alpha$ for $\alpha$-semiconvex functions, the constraints in \eqref{eq:opt-finDim} become bilinear instead of affine.

	\subsection{Least-action kernels with nonnegative Lagrangian and application to inverse optimal control}\label{sec:least-action}
	We assume from now on that $\X$ is a subset of the spacetime vector space $\R\times\R^d$, i.e.\ $x=(t,r)$. Given a Lagrangian function $L:\R\times \R^d\times \R^d\rightarrow\pR$, the action $J\left((t_0,r_0),(t_1,r_1),r(\cdot)\right)$ along an absolutely continuous trajectory %
	${r(\cdot):[t_0,t_1]\rightarrow\R^d}$ going from  $(t_0,r_0)$ to $(t_1,r_1)$ is defined as follows:
	\begin{equation}\label{eq:action_along_traj}
		J\left((t_0,r_0),(t_1,r_1),r(\cdot)\right):=\int_{t_0}^{t_1}L(s,r(s),\dot{r}(s))ds \text{ with } r(t_0)=r_0, \; r(t_1)=r_1.
	\end{equation}
	\begin{Definition}[Maupertuis kernel]	\label{def:Maup_kernel}
		Given a set of absolutely continuous trajectories $S_{\text{traj}}\subset C^0(\R,\R^d)$, we define the Maupertuis kernel $b_{\text{Maup}}:\X\times \X\rightarrow\pmR$, and its asymmetrical version $b^{asym}_{\text{Maup}}$, between two points $x_0=(t_0,r_0)$ and $x_1=(t_1,r_1)$ belonging to $\X$ as follows:
		\begin{align}	\label{eq:Maup_kernel}
			b_{\text{Maup}}(x_0,x_1)&:=-\operatorname{sign}(t_1-t_0)\inf_{\substack{r(\cdot)\in S_{\text{traj}}\\ r(t_0)=r_0\\r(t_1)=r_1}} \underbrace{\int_{t_0}^{t_1}L(s,r(s),\dot{r}(s)) ds,}_{J\left((t_0,r_0),(t_1,r_1),r(\cdot)\right)}\\
			b^{\text{asym}}_{\text{Maup}}(x_0,x_1)&:=(1+\delta^\top_{t_1\ge t_0})b_{\text{Maup}}(x_0,x_1), \label{eq:Maup_kernel_asym}
		\end{align}
		with $\delta^\top_{t_1\ge t_0}=0$ if $t_1\ge t_0$ and $+\infty$ otherwise.
	\end{Definition}
	\begin{Lemma}[Tropical positivity for nonnegative Lagrangian]\label{lem:tpsd_Maup}
		If $L(\cdot)\geq 0$, then $b_{\text{Maup}}$ is tpsd.\footnote{If $S_{\text{traj}}=W^{1,1}(\R,\R^d)$ and $L(s,r,v)=L(v)$, then this statement is even ``if and only if'' using \eqref{eq:Hopf}.}
	\end{Lemma}
	\begin{proof}
		As $L(\cdot)\geq 0$, $b_{\text{Maup}}(\cdot,\cdot)\le 0$, and, for any $x\in\X$, $b_{\text{Maup}}(x,x)=0$. Moreover, for any $x_0,x_1\in\X$, $b_{\text{Maup}}(x_1,x_0)$ is still the infimum over the trajectories joining $x_0$ and $x_1$, and the term $\operatorname{sign}(t_1-t_0)$ compensates the permutation of the bounds in the integral. Hence $b_{\text{Maup}}$ is also symmetric, and we have shown that it is tpsd. 
	\end{proof}
	The evaluation of a Maupertuis kernel amounts to solving an optimal control problem. An approach of optimal control in terms of ``tropical'' kernels has been developed in \citet{KM} and \cite{McEneaney2006}.
	The same kernels have also been studied in the setting
	of ``Lax-Oleinik'' semigroups~\citep{fathi}.
	In general, computing $b_{\text{Maup}}$ is as hard as solving a HJB PDE. However there are cases where the value of $b_{\text{Maup}}$ is known, for instance through Lax-Hopf formulas as recalled in \citet[Theorem 1.3.1]{Cannarsa2004}. This correspond to the case where $L(s,r,v)=L(v)$, $L$ is convex, and $S_{\text{traj}}=W^{1,1}(\R,\R^d)$, in which case \eqref{eq:Maup_kernel} writes simply as a perspective function for $t_1\neq t_0$
	\begin{equation}	\label{eq:Hopf}
		b^{Hopf}_{\text{Maup}}(x_0,x_1):=-|t_1-t_0| L\left(\frac{r_1-r_0}{t_1-t_0}\right)
	\end{equation}
	and $0$ if $x_0=x_1$, $-\infty$ if $r_1\neq r_0$ and $t_0=t_1$. Another special case where $b_{\text{Maup}}$ is known is linear-quadratic optimal control where the kernel $b_{\text{Maup}}$ can be computed by solving differential Riccati equations for time-invariant \citep{dower2015zhang} or time-varying systems \citep{aubin2020hard_control}. 
	Unlike the off-the-shelf Hilbertian kernels used typically in machine learning, we underline that the kernel $b_{\text{Maup}}$ is canonically defined by the triplet $(\X,S_{\text{traj}}, L)$.
	
	Note that requiring as in \Cref{lem:tpsd_Maup} the Lagrangian to take nonnegative values is relevant in applications to optimal control or in the setting of Finsler metrics, in which a distance is defined as the minimum of a positive ``length'' over a set of paths.
	Besides, we can detail the consequence of \Cref{lem:tpsd_Maup} in light of our previous results. The fact that $b_{\text{Maup}}$ vanishes on the diagonal corresponds to the absence of self-loops. Such loops, seen as rewards specific to time points, constitute a translation of the value function, see \eqref{eq:opt-stop-time}-\eqref{eq:value_function} below, and can be incorporated adding a function $\phi$ to $b_{\text{Maup}}$ following \Cref{lem:tpsd_negative-valued}. Moreover, as a consequence of the Bellman principle of optimality, the kernel $b_{\text{Maup}}$ is idempotent (and so is $b^{\text{asym}}_{\text{Maup}}$), whence self-factorized in \Cref{prop:tpsd_feat_map}.
	
	We can identify $\Rg(B^{\text{asym}}_{\text{Maup}})$ based on \eqref{eq:RgB_tpsd}. For every function $f\in\Rg(b^{\text{asym}}_{\text{Maup}})$, there exists $w:\X\mapsto \pR$ such that
	\begin{align}	
		\forall \, x_0=(t_0,r_0)\in\X,\; f(x_0)&=\sup_{x_1\in\X} b^{\text{asym}}_{\text{Maup}}(x_0,x_1) \dotbminus w(x_1)\nonumber \\
		&\hspace{-2cm}=-\inf_{\substack{(t_1,r_1)\in\X,\, t_1\ge t_0,\\ r(\cdot)\in S_{\text{traj}},\, r(t_0)=r_0\\r(t_1)=r_1}} \int_{t_0}^{t_1}L(s,r(s),\dot{r}(s)) ds \dotplus w(t_1,r_1).\label{eq:opt-stop-time}
	\end{align}
	Equation \eqref{eq:opt-stop-time} is precisely, up to a sign, the definition of an optimal stopping time problem, where $w(t_1,r_1)$ is a final cost obtained for choosing to leave the game at $(t_1,r_1)$ \citep{Bensoussan1982,Barles1987}. The player has then to determine, given $w$ and starting from $(t_0,r_0)$, the corresponding final and optimal $(t_1,r_1)$. A formula similar to \eqref{eq:opt-stop-time} can be obtained for $B_{\text{Maup}}$, removing the constraint $t_1\ge t_0$, and thus breaking causality. Nevertheless $b_{\text{Maup}}$ has the advantage of a direct interpretation as a cost for going from one spacetime point to another.
	
	\begin{Example}[Inverse optimal control problem with stopping time]\label{example:stopping_times}
		We now suppose that the terminal cost $w$ in the optimal stopping problem~\eqref{eq:opt-stop-time}
		is unknown, whereas the Lagrangian $L$ is known. Our aim is to infer $w$, and even in fact the value function everywhere, assuming that the value function is only observed at certain points, i.e., that a collection of approximate measurements of the value function, $\bar{y}_m\simeq f(x_m)$, $x_m=(t_m,r_m)\in\X$ with $m\in \mathcal{I}$, is available.   This inverse problem can be cast in
		the form of the representer theorem (\Cref{cor:representer-thm}),
		by considering a loss function %
		\[
		\cL((f(x_m))_{m\in\cI}) = N( (\bar{y}_m-f(x_m))_{m\in\mathcal{I}})
		\]
		where $N$ is any norm.
		Then, by \Cref{cor:representer-thm}, reconstructing the unknown stopping cost  amounts to finding solutions $(p_m,y_m)_{m\in \cI}\in (\X\times\R)^\cI$ minimizing $\cL((y_m)_{m\in\cI})$ under the constraints
		\begin{gather*}
			\forall \, n,m\in\cI,\, y_n-y_m\ge b^{\text{asym}}_{\text{Maup}}(x_n,p_m)\dotbminus b^{\text{asym}}_{\text{Maup}}(x_m,p_m).%
		\end{gather*}
		Using that $b_{\text{Maup}}$ is idempotent, we have that $b^{\text{asym}}_{\text{Maup}}(x_n,p_m)\dotbminus b^{\text{asym}}_{\text{Maup}}(x_m,p_m)\ge b^{\text{asym}}_{\text{Maup}}(x_n,x_m)$, achieved for $p_m=x_m$, whence the problem reduces to minimizing the function 
		$\cL$ over $(y_m)_{m\in \cI}\in \R^\cI$ such that
		\begin{gather*}
			\forall \, n,m\in\cI,\, y_n-y_m\ge b^{\text{asym}}_{\text{Maup}}(x_n,x_m).%
		\end{gather*}		\Cref{prop:finite-solution}-\ref{it_finiteinterp} then provides an expression of the corresponding value function. Indeed, if $(y^*_m)_{m\in\cI}$ is an optimal solution of the above problem, an admissible stopping cost is simply $w=\inf_{m\in\cI}[\delta^\top_{x_m}(\cdot)-y^*_m]$, in other words, $w(x_m)=-y^*_m$ for all $m\in\cI$ and $w$ is $+\infty$ elsewhere. Many variations of this problem can be considered, introducing other loss functions, e.g. to infer bounds on the value function at a certain point while knowing measurements of the same value function at other points.
	\end{Example}%
	
	Another interesting case is when the terminal time $T$ is given. As a matter of fact, given a terminal cost $\psi_T(\cdot)$, setting $\X=(-\infty,T]\times \R^d$, recall that the value function is defined as
	\begin{gather}\label{eq:value_function}
		\bar{V}(t_0,r_0)=\inf_{\substack{r_T\in\R^d,\,  r(\cdot)\in S_{\text{traj}}\\ r(t_0)=r_0,\, r(T)=r_T}} \int_{t_0}^{T}L(s,r(s),\dot{r}(s)) ds+\psi_T(r_T).
	\end{gather}
	This formula can actually be rewritten as
	\begin{gather}\label{eq:HJB_variationnel}
		\bar{V}(\cdot,\cdot)=-\min\{ \text{$\hat{V}(\cdot,\cdot)\in \Rg(B_{\text{Maup}})$ s.t.\ $\hat{V}(T,r)=\psi_T(r)$ for all $r$}\}.
	\end{gather}%
	Indeed, a solution is merely $\bar{V}(\cdot,\cdot)=\bB_{\text{Maup}}\tilde{\psi}$ where $\tilde{\psi}(t,r)=\psi_T(r)+\delta^\top_T(t)$ as, for any $\hat{V}(\cdot,\cdot)$ in the argmin of \eqref{eq:HJB_variationnel}, since $b_{\text{Maup}}$ is idempotent,
	\begin{align*}
		\hat{V}(t_0,r_0)&=\sup_{(t,r)\in \R^{d+1}} b_{\text{Maup}}((t_0,r_0),(t,r))-\hat{V}(t,r)\\
		&\stackrel{t=T}{\ge}\sup_{r\in\R^d} b_{\text{Maup}}((t_0,r_0),(T,r))-\psi_T(r)=-\bar{V}(t_0,r_0)
	\end{align*}
	This recovers the well-known property that the value function is the largest subsolution of the HJB equation. \\

	\noindent \tb{Space-restrictions of $b_{\text{Maup}}$}: With fixed $t_0,t_1\in\R$, one can also consider a kernel $b^{[t_0,t_1]}_{\text{Maup}}(r_0,r_1):=b_{\text{Maup}}(x_0,x_1)$. This kernel defined on $\R^d\times\R^d$ is symmetric only in specific cases, e.g.\ a locally reversible set of trajectories $S_{\text{traj}}$ and a time-independent Lagrangian such that $L(r,v)= L(r,-v)$ for all $(r,v)$. 
	These kernels, which are the most commonly used in the ``tropical approach'' of optimal control \citep{KM,McEneaney2006,dower2015zhang}, are 
	not idempotent. 
	\begin{Example}[Inverse optimal control problem: interpolation of the value function with a fixed final time]\label{example:final_time} We now assume that an initial time $t_0$ and a final time $T$ are fixed, and that
		the Lagrangian $L$ is known. Our aim is to recover
		an unknown terminal cost $\psi_T$, and thus, the value function $\bar{V}$ everywhere, given measurements $(\bar{y}_m)_{m\in \cI}\in \R^\cI$ 
		at sample points. For instance, if we consider for simplicity the exact interpolation problem, we may assume that $\bar{y}_m=-\bar{V}(t_0,r_m)$, at some point $r_m$
		(the variant with approximate measurements can be handled as per~\Cref{example:stopping_times}). %
		One can apply \Cref{prop:finite-solution}, so the interpolation problem corresponds to finding solutions $(p_m)_{m\in \cI}\in (\R^d)^\cI$  of  
		\begin{gather*}
			\forall\,  n,m\in\cI,\, \bar{y}_n-\bar{y}_m\ge b^{[t_0,T]}_{\text{Maup}}(r_n,p_m)\dotbminus b^{[t_0,T]}_{\text{Maup}}(r_m,p_m).
		\end{gather*}
		Here each $p_m$ can be interpreted as a point to reach at time $T$ starting from $(t_0,r_m)$. Given such $p_m$, an admissible terminal cost can be reconstructed as $\psi_T(r)=\inf_{m\in\cI}[\delta^\top_{p_m}(r)\dotplus b(x_m,p_m)-\bar y_m]$.%
	\end{Example}

	\bibliographystyle{apalike}
	\bibliography{Kernels_TropGeometry}	

\begin{thebibliography}{}

\bibitem[Akian, 1999]{akiantams}
Akian, M. (1999).
\newblock Densities of idempotent measures and large deviations.
\newblock {\em Transactions of the American Mathematical Society},
  351(11):4515--4543.

\bibitem[Akian et~al., 2005]{Akian04setcoverings}
Akian, M., Gaubert, S., and Kolokoltsov, V.~N. (2005).
\newblock Set coverings and invertibility of functional {Galois} connections.
\newblock In Litvinov, G.~L. and Maslov, V.~P., editors, {\em Idempotent
  Mathematics and Mathematical Physics}, Contemporary Mathematics, pages
  19--51. American Mathematical Society.

\bibitem[Akian et~al., 2008]{asma-art-1}
Akian, M., Gaubert, S., and Lakhoua, A. (2008).
\newblock The max-plus finite element method for solving deterministic optimal
  control problems: basic properties and convergence analysis.
\newblock {\em SIAM J. Control Optim.}, 47(2):817--848.

\bibitem[Akian et~al., 1994]{bellman}
Akian, M., Quadrat, J., and Viot, M. (1994).
\newblock Bellman processes.
\newblock In {\em 11th International Conference on Analysis and Optimization of
  Systems : Discrete Event Systems}, volume 199 of {\em Lecture notes in
  control and information sciences}. Springer Verlag.

\bibitem[Aronszajn, 1950]{aronszajn50theory}
Aronszajn, N. (1950).
\newblock Theory of reproducing kernels.
\newblock {\em Transactions of the American Mathematical Society}, 68:337--404.

\bibitem[Artstein-Avidan and Milman, 2009]{ArtsteinAvidan2009}
Artstein-Avidan, S. and Milman, V. (2009).
\newblock The concept of duality in convex analysis, and the characterization
  of the legendre transform.
\newblock {\em Annals of Mathematics}, 169(2):661--674.

\bibitem[Aubin-Frankowski, 2021a]{aubin2020Riccati}
Aubin-Frankowski, P.-C. (2021a).
\newblock Interpreting the dual {Riccati} equation through the {LQ} reproducing
  kernel.
\newblock {\em Comptes Rendus. Math{\'{e}}matique}, 359(2):199--204.

\bibitem[Aubin-Frankowski, 2021b]{aubin2020hard_control}
Aubin-Frankowski, P.-C. (2021b).
\newblock Linearly constrained linear quadratic regulator from the viewpoint of
  kernel methods.
\newblock {\em {SIAM} Journal on Control and Optimization}, 59(4):2693--2716.

\bibitem[Baccelli et~al., 1992]{bcoq}
Baccelli, F., Cohen, G., Olsder, G., and Quadrat, J. (1992).
\newblock {\em Synchronization and Linearity}.
\newblock Wiley.

\bibitem[Barles and Perthame, 1987]{Barles1987}
Barles, G. and Perthame, B. (1987).
\newblock Discontinuous solutions of deterministic optimal stopping time
  problems.
\newblock {\em {ESAIM}: Mathematical Modelling and Numerical Analysis},
  21(4):557--579.

\bibitem[Bensoussan and Lions, 1982]{Bensoussan1982}
Bensoussan, A. and Lions, J.-L. (1982).
\newblock {\em Applications of variational inequalities in stochastic control}.
\newblock Studies in mathematics and its applications. Elsevier Science.

\bibitem[Berg et~al., 1984]{Berg1984}
Berg, C., Christensen, J. P.~R., and Ressel, P. (1984).
\newblock {\em Harmonic Analysis on Semigroups}.
\newblock Springer New York.

\bibitem[Boyd and Vandenberghe, 2004]{boyd_vandenberghe_2004}
Boyd, S. and Vandenberghe, L. (2004).
\newblock {\em Convex Optimization}.
\newblock Cambridge University Press.

\bibitem[Burkard et~al., 1996]{Burkard1996}
Burkard, R.~E., Klinz, B., and Rudolf, R. (1996).
\newblock Perspectives of {Monge} properties in optimization.
\newblock {\em Discrete Applied Mathematics}, 70(2):95--161.

\bibitem[Butkovi{\v{c}} et~al., 2007]{BSS}
Butkovi{\v{c}}, P., Schneider, H., and Sergeev, S. (2007).
\newblock Generators, extremals and bases of max cones.
\newblock {\em Linear Algebra Appl.}, 421(2-3):394--406.

\bibitem[Calafiore et~al., 2020]{1905.08503}
Calafiore, G.~C., Gaubert, S., and Possieri, C. (2020).
\newblock A universal approximation result for difference of log-sum-exp neural
  networks.
\newblock {\em IEEE Trans. Neural Networks Learn. Syst.}, 31(12):5603--5612.

\bibitem[Cannarsa and Sinestrari, 2004]{Cannarsa2004}
Cannarsa, P. and Sinestrari, C. (2004).
\newblock {\em Semiconcave Functions, Hamilton{\textemdash}Jacobi Equations,
  and Optimal Control}.
\newblock Birkh\"{a}user Boston.

\bibitem[Cartwright and Chan, 2012]{Cartwright2012}
Cartwright, D. and Chan, M. (2012).
\newblock Three notions of tropical rank for symmetric matrices.
\newblock {\em Combinatorica}, 32(1):55--84.

\bibitem[Chancelier and De~Lara, 2021]{chancelier2021capra}
Chancelier, J.-P. and De~Lara, M. (2021).
\newblock Capra-convexity, convex factorization and variational formulations
  for the $\ell_0$ pseudonorm.
\newblock {\em Set-Valued and Variational Analysis}.
\newblock on line.

\bibitem[Cohen et~al., 2004]{Cohen2004}
Cohen, G., Gaubert, S., and Quadrat, J.-P. (2004).
\newblock Duality and separation theorems in idempotent semimodules.
\newblock {\em Linear Algebra and its Applications}, 379:395--422.

\bibitem[Develin and Sturmfels, 2004]{DS04}
Develin, M. and Sturmfels, B. (2004).
\newblock Tropical convexity.
\newblock {\em Doc. Math.}, 9:1--27.
\newblock (Erratum pp.~205--206).

\bibitem[Di~Marino et~al., 2017]{di2017optimal}
Di~Marino, S., Gerolin, A., and Nenna, L. (2017).
\newblock Optimal transportation theory with repulsive costs.
\newblock {\em Topological optimization and optimal transport}, 17:204--256.

\bibitem[Dower and McEneaney, 2015]{Dower2015}
Dower, P.~M. and McEneaney, W.~M. (2015).
\newblock A max-plus dual space fundamental solution for a class of operator
  differential {Riccati} equations.
\newblock {\em {SIAM} Journal on Control and Optimization}, 53(2):969--1002.

\bibitem[Dower and Zhang, 2015]{dower2015zhang}
Dower, P.~M. and Zhang, H. (2015).
\newblock A new fundamental solution for differential riccati equations arising
  in l2-gain analysis.
\newblock In {\em 2015 5th Australian Control Conference (AUCC)}, pages 65--68.

\bibitem[Fathi, 2008]{fathi}
Fathi, A. (2008).
\newblock The weak-{KAM} theorem in lagrangian dynamics.
\newblock Version 10, available from
  \url{https://www.math.u-bordeaux.fr/~pthieull/Recherche/KamFaible/Publications/Fathi2008_01.pdf}.

\bibitem[Feydy et~al., 2019]{feydy2019interpolating}
Feydy, J., S\'{e}journ\'{e}, T., Vialard, F.-X., Amari, S.-i., Trouve, A., and
  Peyr\'{e}, G. (2019).
\newblock Interpolating between optimal transport and {MMD} using {Sinkhorn}
  divergences.
\newblock In {\em International Conference on Artificial Intelligence and
  Statistics (AISTATS)}, volume~89, pages 2681--2690.

\bibitem[Fleming and McEneaney, 2000]{a5}
Fleming, W.~H. and McEneaney, W.~M. (2000).
\newblock A max-plus-based algorithm for a {H}amilton-{J}acobi-{B}ellman
  equation of nonlinear filtering.
\newblock {\em SIAM J. Control Optim.}, 38(3):683--710.

\bibitem[Gaubert and Katz, 2007]{GK}
Gaubert, S. and Katz, R. (2007).
\newblock The {M}inkowski theorem for max-plus convex sets.
\newblock {\em Linear Algebra and Appl.}, 421:356--369.

\bibitem[Gaubert and Niv, 2018]{adiniv}
Gaubert, S. and Niv, A. (2018).
\newblock Tropical totally positive matrices.
\newblock {\em Journal of Algebra}, 515:511 -- 544.

\bibitem[Gierz et~al., 2003]{continuous}
Gierz, G., Hofmann, K., Keimel, K., Lawson, J., Mislove, M., and Scott, D.
  (2003).
\newblock {\em Continuous Lattices and Domains}.
\newblock Encyclopedia of Mathematics and its Applications. Cambridge
  University Press.

\bibitem[Kimeldorf and Wahba, 1971]{kimeldorf1971tcheb}
Kimeldorf, G. and Wahba, G. (1971).
\newblock Some results on {Tchebycheffian} spline functions.
\newblock {\em Journal of Mathematical Analysis and Applications},
  33(1):82--95.

\bibitem[Kolokoltsov and Maslov, 1997]{KM}
Kolokoltsov, V.~N. and Maslov, V.~P. (1997).
\newblock {\em Idempotent analysis and its applications}, volume 401 of {\em
  Mathematics and its Applications}.
\newblock Kluwer Academic Publishers Group.

\bibitem[Lam and Postnikov, 2007]{AlcovedPolytopes}
Lam, T. and Postnikov, A. (2007).
\newblock Alcoved polytopes. {I}.
\newblock {\em Discrete Comput. Geom.}, 38(3):453--478.

\bibitem[Litvinov, 2011]{litvinov2011trop}
Litvinov, G. (2011).
\newblock Tropical mathematics, idempotent analysis, classical mechanics and
  geometry.
\newblock In {\em Spectral Theory and Geometric Analysis}, pages 159--186.
  American Mathematical Society.

\bibitem[Litvinov, 2005]{litvinov2005maslov}
Litvinov, G.~L. (2005).
\newblock Maslov dequantization, idempotent and tropical mathematics: A brief
  introduction.
\newblock {\em Journal of Mathematical Sciences}, 140:426--444.

\bibitem[Maragos et~al., 2021]{maragos2021}
Maragos, P., Charisopoulos, V., and Theodosis, E. (2021).
\newblock Tropical geometry and machine learning.
\newblock {\em Proceedings of the {IEEE}}, 109(5):728--755.

\bibitem[Martinez-Legaz and Singer, 1990]{martinezlegaz1990dualities}
Martinez-Legaz, J. and Singer, I. (1990).
\newblock Dualities between complete lattices.
\newblock {\em Optimization}, 21(4):481--508.

\bibitem[Martinez-Legaz and Singer, 1995]{MartinezLegaz1995subdifferentials}
Martinez-Legaz, J.-E. and Singer, I. (1995).
\newblock Subdifferentials with respect to dualities.
\newblock {\em {ZOR} Zeitschrift für Operations Research Methods and Models of
  Operations Research}, 42(1):109--125.

\bibitem[Mary, 2005]{mary2005theory}
Mary, X. (2005).
\newblock Theory of subdualities.
\newblock {\em Journal d'Analyse Math{\'{e}}matique}, 97(1):203--241.

\bibitem[McEneaney, 2006]{McEneaney2006}
McEneaney, W.~M. (2006).
\newblock {\em Max-Plus Methods for Nonlinear Control and Estimation}.
\newblock Birkh\"{a}user-Verlag.

\bibitem[McEneaney, 2007]{curseofdim}
McEneaney, W.~M. (2007).
\newblock A curse-of-dimensionality-free numerical method for solution of
  certain {HJB} {PDE}s.
\newblock {\em SIAM J. Control Optim.}, 46(4):1239--1276.

\bibitem[Mont{\'{u}}far et~al., 2022]{2104.08135}
Mont{\'{u}}far, G., Ren, Y., and Zhang, L. (2022).
\newblock Sharp bounds for the number of regions of maxout networks and
  vertices of minkowski sums.
\newblock {\em {SIAM} Journal on Applied Algebra and Geometry}, 6(4):618--649.

\bibitem[Moreau, 1970]{moreau1970infconvol}
Moreau, J.~J. (1970).
\newblock {Inf-convolution, sous-additivit{\'e}, convexit{\'e} des fonctions
  num{\'e}riques}.
\newblock {\em {Journal de Math{\'e}matiques Pures et Appliqu{\'e}es}}, pages
  109--154.

\bibitem[Ong et~al., 2004]{ong2004learning}
Ong, C.~S., Mary, X., Canu, S., and Smola, A.~J. (2004).
\newblock Learning with non-positive kernels.
\newblock In {\em International Conference on Machine Learning (ICML)}, pages
  81--88. {ACM} Press.

\bibitem[Pallaschke and Rolewicz, 1997]{Pallaschke1997}
Pallaschke, D. and Rolewicz, S. (1997).
\newblock {\em Foundations of Mathematical Optimization}.
\newblock Springer Netherlands.

\bibitem[Rachev and Rüschendorf, 1998]{Rachev1998}
Rachev, S.~T. and Rüschendorf, L. (1998).
\newblock {\em Mass Transportation Problems}.
\newblock Springer-Verlag.

\bibitem[Saitoh and Sawano, 2016]{saitoh16theory}
Saitoh, S. and Sawano, Y. (2016).
\newblock {\em Theory of Reproducing Kernels and Applications}.
\newblock Springer Singapore.

\bibitem[Santambrogio, 2017]{santambrogio2017}
Santambrogio, F. (2017).
\newblock $\lbrace${E}uclidean, metric, and {W}asserstein$\rbrace$ gradient
  flows: an overview.
\newblock {\em Bulletin of Mathematical Sciences}, 7(1):87--154.

\bibitem[Schölkopf et~al., 2001]{scholkopf2001representer}
Schölkopf, B., Herbrich, R., and Smola, A.~J. (2001).
\newblock A generalized representer theorem.
\newblock In {\em Computational Learning Theory (COLT)}, pages 416--426.
  Springer Berlin Heidelberg.

\bibitem[Sch{\"o}lkopf and Smola, 2002]{scholkopf02learning}
Sch{\"o}lkopf, B. and Smola, A. (2002).
\newblock {\em Learning with Kernels: Support Vector Machines, Regularization,
  Optimization, and Beyond}.
\newblock MIT Press.

\bibitem[Schwartz, 1964]{schwartz1964sev}
Schwartz, L. (1964).
\newblock Sous-espaces hilbertiens d'espaces vectoriels topologiques et noyaux
  associ\'{e}s (noyaux reproduisants).
\newblock {\em Journal d’Analyse Mathématique}, 13:115--256.

\bibitem[Seijo and Sen, 2011]{Seijo2011}
Seijo, E. and Sen, B. (2011).
\newblock Nonparametric least squares estimation of a multivariate convex
  regression function.
\newblock {\em The Annals of Statistics}, 39(3):1633--1657.

\bibitem[Singer, 1984]{singer1984conj}
Singer, I. (1984).
\newblock Conjugation operators.
\newblock In {\em Lecture Notes in Economics and Mathematical Systems}, pages
  80--97. Springer Berlin Heidelberg.

\bibitem[Singer, 1997]{singer1997}
Singer, I. (1997).
\newblock {\em Abstract Convex Analysis}.
\newblock Wiley-Interscience and Canadian Mathematics Series of Monographs and
  Texts. Wiley-Interscience, 1 edition.

\bibitem[Steinwart and Christmann, 2008]{steinwart08support}
Steinwart, I. and Christmann, A. (2008).
\newblock {\em Support Vector Machines}.
\newblock Springer.

\bibitem[Taylor et~al., 2016]{Taylor2016}
Taylor, A.~B., Hendrickx, J.~M., and Glineur, F. (2016).
\newblock Smooth strongly convex interpolation and exact worst-case performance
  of first-order methods.
\newblock {\em Mathematical Programming}, 161(1-2):307--345.

\bibitem[Tran, 2020]{Tran2020}
Tran, N.~M. (2020).
\newblock Tropical {Gaussians}: a brief survey.
\newblock {\em Algebraic Statistics}, 11(2):155--168.

\bibitem[Villani, 2003]{Villani2003}
Villani, C. (2003).
\newblock {\em Topics in optimal transportation}.
\newblock American Mathematical Society.

\bibitem[Villani, 2009]{Villani2009}
Villani, C. (2009).
\newblock {\em Optimal Transport}.
\newblock Springer Berlin Heidelberg.

\bibitem[Volle et~al., 2013]{volle2013}
Volle, M., Mart{\'{\i}}nez-Legaz, J.~E., and Vicente-P{\'{e}}rez, J. (2013).
\newblock Duality for closed convex functions and evenly convex functions.
\newblock {\em Journal of Optimization Theory and Applications},
  167(3):985--997.

\bibitem[Wei{\ss} et~al., 2016]{Weiss2016}
Wei{\ss}, C., Knust, S., Shakhlevich, N., and Waldherr, S. (2016).
\newblock The assignment problem with nearly {Monge} arrays and incompatible
  partner indices.
\newblock {\em Discrete Applied Mathematics}, 211:183--203.

\bibitem[Yoshida et~al., 2023]{yoshida2021TropicalSV}
Yoshida, R., Takamori, M., Matsumoto, H., and Miura, K. (2023).
\newblock Tropical support vector machines: Evaluations and extension to
  function spaces.
\newblock {\em Neural Networks}, 157:77--89.

\bibitem[Yoshida et~al., 2019]{yoshidatropPCAphylo}
Yoshida, R., Zhang, L., and Zhang, X. (2019).
\newblock Tropical principal component analysis and its application to
  phylogenetics.
\newblock {\em Bulletin of Mathematical Biology}, 81:568--597.

\bibitem[Yu, 2014]{yu2014tropicalizing}
Yu, J. (2014).
\newblock Tropicalizing the positive semidefinite cone.
\newblock {\em Proceedings of the American Mathematical Society},
  143(5):1891--1895.

\bibitem[Zhang et~al., 2018]{pmlr-v80-zhang18i}
Zhang, L., Naitzat, G., and Lim, L.-H. (2018).
\newblock Tropical geometry of deep neural networks.
\newblock In {\em International Conference on Machine Learning (ICML)},
  volume~80, pages 5824--5832.

\end{thebibliography}
\end{document}